%% file: paper.tex
\documentclass[%
leqno,%
10pt%
]{article}

\usepackage[utf8]{inputenc}
\usepackage[T1]{fontenc}
\usepackage[english]{babel}

\usepackage[intlimits]{amsmath}
\usepackage{amssymb, amsfonts}
\usepackage[amsmath, thmmarks]{ntheorem}

\usepackage[intlimits]{amsmath}
\usepackage{amssymb, amsfonts}
\usepackage[amsmath]{ntheorem}

\usepackage{bbm}

\usepackage[all]{xy}

\numberwithin{equation}{section}

\newtheorem{theoremcounter}{theoremcounter}[section]

\theoremstyle{plain}
\theoremseparator{.}
\newtheorem{remark}{Remark}

\theoremstyle{plain}
\newtheorem{definition}[theoremcounter]{Definition}
\newtheorem{lemma}[theoremcounter]{Lemma}
\newtheorem{theorem}[theoremcounter]{Theorem}
\newtheorem{corollary}[theoremcounter]{Corollary}
\newtheorem{proposition}[theoremcounter]{Proposition}

\theorembodyfont{\normalfont}
\theoremsymbol{\ensuremath{\Box}}
\newtheorem*{proof}{Proof}

\newcommand{\ZZ}{\ensuremath{\mathbb{Z}}}
\newcommand{\QQ}{\ensuremath{\mathbb{Q}}}
\newcommand{\RR}{\ensuremath{\mathbb{R}}}
\newcommand{\CC}{\ensuremath{\mathbb{C}}}

\newcommand{\Hom}{\ensuremath{\mathrm{Hom}}}

\newcommand{\Mat}[2]{\ensuremath{\mathrm{M}_{#1}(#2)}}
\newcommand{\MatT}[2]{\ensuremath{\mathrm{M}^\T_{#1}(#2)}}

\newcommand{\GL}[1]{\ensuremath{\mathrm{GL}_{#1}}}
\newcommand{\SL}[1]{\ensuremath{\mathrm{SL}_{#1}}}
\newcommand{\Sp}[1]{\ensuremath{\mathrm{Sp}_{#1}}}

\newcommand{\T}{\ensuremath{\mathrm{T}}}

\newcommand{\tr}{\ensuremath{\mathrm{tr}}}

\newcommand{\slashdiv}{\ensuremath{\mathop{/}}}

\renewcommand{\Re}{\ensuremath{\mathfrak{Re}}}
\renewcommand{\Im}{\ensuremath{\mathfrak{Im}}}

\newcommand{\HS}{\mathbb{H}}

\newcommand{\wtM}{\ensuremath{\mathcal{M}}}
\newcommand{\wtW}{\ensuremath{\mathcal{W}}}

\def\frak{\mathfrak}

\def\dg{{\frak g}}
\def\dgl{\dg\dl}
\def\dh{{\frak h}}

\def\dk{{\frak k}}
\def\dl{{\frak l}}
\def\dm{{\frak m}}

\def\dgo{{\frak o}}

\def\dr{{\frak r}}
\def\ds{{\frak s}}
\def\dsl{\ds\dl}

\def\dU{{\frak U}}

\def\dz{{\frak z}}

\def\dZ{{\frak Z}}

\def\cal{\mathcal}

\def\cC{{\cal C}} 

\def\Bbb{\mathbb}

\def\bC{\Bbb C}
\def\bD{\Bbb D}

\def\bJ{\Bbb J}

\def\bR{\Bbb R}

\def\rmM{\mathop{\rm M}\nolimits}

\def\sL{\mbox{\it \L}}

\newcommand{\Jac}[1]{\ensuremath{G^J_{#1}}}
\newcommand{\JacF}[1]{\ensuremath{\Gamma^J_{#1}}}
\newcommand{\JacFpara}[1]{\ensuremath{\Gamma^J_{#1,\infty}}}

\newcommand{\tJac}[1]{\ensuremath{\t G^J_{#1}}}
\newcommand{\tJacK}[1]{\ensuremath{\t K^J_{#1}}}

\def\oh{{\ts\frac{1}{2}}}
\def\ov{\overline}
\def\pz{\partial_z}
\def\pbz{\partial_{\ov z}}
\def\scrm{\scriptsize\rm}
\def\SO{\mathop{\rm SO}\nolimits}

\def\Sp{\mathop{\rm Sp}\nolimits}
\def\Span{\mathop{\rm Span}\nolimits}
\def\t{\tilde}
\def\thup{{\mbox{\scrm th}}}

\def\ts{\textstyle}

\def\dog{differential operator}
\def\ido{invariant \dog}
\def\iff{if and only if}
\def\rep{representation}
\def\uea{universal enveloping algebra}
\begin{document}
\title{Harmonic Maa\ss-Jacobi forms of degree 1 with higher rank indices}
\author{
\begin{minipage}{13.5em} \centering{Charles Conley}\\
\small{conley@unt.edu} \vspace{1.5em}\\
Department of Mathematics \\ University of North Texas \\ Denton TX 76203-1430 USA
\end{minipage}
\ and\ 
\begin{minipage}{13.5em} \centering{Martin Raum}\\
\small{mraum@mpim-bonn.mpg.de} \vspace{1.5em}\\
MPI f\"ur Mathematik \\ Vivatsgasse 7 \\ 53111 Bonn, Germany
\end{minipage}
}

\date{\today}

\maketitle

\input{abstract}
\input{introduction}
\input{maassjacobiforms}
\input{mockjacobiforms}
\input{semiholomorphicforms}
\input{casimiroperators}

\bibliographystyle{amsalpha}
\bibliography{bibliography}

\end{document}

%% file: abstract.tex
\begin{abstract}
We define and investigate real analytic weak Jacobi forms of degree~1 and arbitrary rank.  En route we calculate the Casimir operator associated to the maximal central extension of the real Jacobi group, which for rank exceeding~1 is of order~4.  In ranks exceeding~1, the notions of H-harmonicity and semi-holomorphicity are the same.

\bigbreak\noindent

\bigbreak\noindent
The first author was partially supported by Simons Foundation Collaboration Grant 207736.  The second author was partially supported by Max Planck Institute for Mathematics.

\end{abstract}

%% file: introduction.tex
\section{Introduction}
\label{sec:introduction}

The theory of holomorphic Jacobi forms was developed by Eichler and Zagier in the course of their work on the Saito-Kurokawa conjecture \cite{EZ85}.  Later Berndt and Schmidt initiated a theory of real analytic Jacobi forms \cite{BS98}, which was developed further by Pitale \cite{Pi09}.  In the real analytic case, holomorphicity is replaced by the requirement that the forms be eigenfunctions of the Casimir operator, a third order operator which generates the center of the algebra of invariant operators \cite{BCR10}.

Bringmann and Richter studied harmonic Maa\ss-Jacobi forms in the sense of Pitale \cite{BR10}, but with a weak growth condition that includes the $\mu$-function discovered by Zwegers.  Zwegers had used this function in \cite{Zw02} to understand the hitherto mysterious mock modular forms discovered by Ramanujan in the early~20$^\textrm{th}$ century.  His work has been the focus of intense interest, having applications to mock theta functions, combinatorics, and physics.

Dabholkar, Murthy, and Zagier discovered deep connections between Jacobi forms and the state-counting properties of quantum black holes \cite{DMZ11}.  Since they begin with meromorphic Siegel modular forms, they need only study a very special class of Maa\ss-Jacobi forms: those which are semi-holomorphic.  Semi-holomorphic forms also appeared as ``Fourier Jacobi coefficients'' of particular Siegel Maass forms \cite{BRR12}.  This illustrates the importance of focusing particular attention on such forms.

Zwegers generalized the $\mu$-function to higher Jacobi forms \cite{Zw10} by demonstrating the modularity of the completed multivariable Appell functions arising from certain character formulas for Lie superalgebras \cite{KW94, KW01, STT05}.  Semi-holomorphic forms do not include the $\mu$-function, so it is necessary to consider the more general of notion of H-harmonicity introduced in \cite{BRRZ11} in the rank~1 case.  It is founded on a Laplace-like operator acting on the elliptic variable of a Maa\ss-Jacobi form.

In the present work we generalize the notion of harmonic weak Maa\ss-Jacobi forms of degree~1 to arbitrary indices of higher rank, in a manner which includes the Appell functions treated in \cite{Zw10}.  Let $\Jac{N}$ be the rank~$N$ Jacobi group $\SL{2}(\RR) \ltimes (\RR^N \otimes \RR^2)$, the semidirect product action being trivial on the first factor, and let $\tJac{N}$ be its central extension by the additive group $\MatT{N}{\RR}$ of real symmetric $N \times N$ matrices.  An important ingredient of our work is the center of the universal enveloping algebra of $\tJac{N}$.  Using ideas developed by Borho \cite{Bo76}, Quesne \cite{Qu88}, and Campoamor-Stursburg and Low \cite{CL09}, we prove in Section~\ref{sec:casimiroperators} that this center is the polynomial algebra generated by $\MatT{N}{\RR}$ and one additional element $\t\Omega_N$ of degree~$N + 2$.  We refer to $\t\Omega_N$ as the {\em Casimir element\/} of $\tJac{N}$, as it is in some sense a lift of the Casimir element of $\SL{2}$.

Given any action of $\tJac{N}$, we refer to the operator by which $\t\Omega_N$ acts as the {\em Casimir operator\/} with respect to the action.  In Section~\ref{sec:maassjacobiforms} we give formulas for the Casimir operators with respect to the standard slash actions of $\tJac{N}$, in terms of both the usual coordinates~(\ref{casop}) and the raising and lowering operators~(\ref{eq:casimir_raisinglowering}).  For $N \ge 2$ these operators are of order~4.

Let $\Gamma^J_N$ be the {\em full Jacobi group,\/} the integer points of $G^J_N$.  The slash actions of $\t G^J_N$ of interest here all drop to actions of $\Gamma^J_N$.  We define a {\em Maa\ss-Jacobi form\/} to be an eigenform of the Casimir operator with respect to such an action, invariant under $\Gamma^J_N$ and satisfying a certain growth condition.

In Section~\ref{sec:mockjacobiforms} we deduce the natural generalization of H-harmonicity to all of $\tJac{N}$, starting from three fundamental requirements.  Surprisingly, it is much stronger than expected: for $N \ge 2$ it coincides with the notion of semi-holomorphicity discussed in the subsequent section.

In Section~\ref{sec:semiholomorphicforms} we investigate a distinguished subspace of the space of Maa\ss-Jacobi forms, the space of {\em semi-holomorphic forms.\/}  We show that in the higher rank case it is connected to the space of skew-holomorphic Jacobi forms: we define a $\xi$-operator~(\ref{eq:definition_xioperator}) which maps any harmonic Maa\ss-Jacobi form to the derivation of its non-holomorphic part, a skew-holomorphic Jacobi form in the sense of \cite{Sk90, Ha06}.  In Corollary~\ref{cor:xisurjectivityandimages} we show that all possible cuspidal non-holomorphic parts occur.

The Zagier-type dualities proved in Corollary~\ref{cor:zagiertypeduality} demonstrate the arithmetic relevance of our more general construction.  As Bringmann and Richter remark in the rank~1 case \cite{BR10}, this relates holomorphic parts not only to one another, but also to non-holomorphic parts.

The paper concludes with Section~\ref{sec:casimiroperators}, in which we use the algorithm developed by Helgason \cite{Helg77} to deduce the invariant and covariant \dog s presented in Section~\ref{sec:maassjacobiforms}.

\subsection*{Acknowledgements}
\textit{The authors thank Olav Richter and Don Zagier for inspiring discussions.}


%% file: maassjacobiforms.tex
\section{Maa\ss-Jacobi forms with lattice indices}
\label{sec:maassjacobiforms}

We first fix some notation.  All vector spaces are complex unless we indicate otherwise.  Let $\Mat{m, n}{R}$ denote the space of $m \times n$ matrices over a ring $R$, abbreviate $\Mat{n, n}{R}$ as $\Mat{n}{R}$, and let $\MatT{n}{R}$ be the symmetric subspace of $\Mat{n}{R}$.  Write $A^\T$ and $\tr(A)$ for the transpose and (when $A$ is square) trace of a matrix $A$, respectively.  Regarding elements of $R^m$ as column vectors, we will freely identify $R^m \otimes R^n$ with $\Mat{m, n}{R}$ via $v \otimes w \mapsto v w^\T$.  

Write $\epsilon_i$ for the $i^\thup$ standard basis vector of $R^m$ and $\epsilon_{ij}$ for the elementary matrix with $(i, j)^\thup$ entry~$1$ and other entries~$0$, the sizes of $\epsilon_i$ and $\epsilon_{ij}$ being determined by context.  Let $I_n$ be the identity matrix in $\Mat{n}{R}$, and set
$$ J_{2n} := \Bigl(\begin{matrix} 0 & -I_n \\ I_n & \ 0\end{matrix} \Bigr). $$

The real Jacobi group $\Jac{N}$ of rank~$N$ and its subgroup $\JacF{N}$, the {\em full Jacobi group,\/} are
\begin{equation} \label{eq:jacandfulljac}
   \Jac{N} := \SL{2}(\RR) \ltimes (\RR^N \otimes \RR^2),
   \qquad
   \JacF{N} := \SL{2}(\ZZ) \ltimes (\ZZ^N \otimes \ZZ^2).
\end{equation}
The product in $\Jac{N}$ arises from the natural right action of $\SL{2}(\RR)$ on $\RR^2$.  It can be written most simply using the above identification of $\RR^N \otimes \RR^2$ with $\Mat{N, 2}{\RR}$: for $M$, $\check M \in \SL{2}(\RR)$ and $X$, $\check X \in \Mat{N, 2}{\RR}$,
$$ (M, X) (\check M, \check X) = (M \check M, X \check M + \check X). $$

Let $\HS := \{\tau \in \CC : \Im(\tau) > 0\}$ be the Poincar\'e upper half plane, and define 
$$ \HS_{1,N} := \HS \times \CC^N. $$
We will write $\tau := x + iy$ for the $\HS$-coordinate and $z_j := u_j + iv_j$ for the $\CC^N$-coordinates.  We will be interested in a certain family of slash actions ({\it i.e.,\/} right  actions) of $\JacF{N}$ on $C^{\infty}(\HS_{1,N})$.  These actions are not restrictions of actions of $\Jac{N}$, but rather quotients of restrictions of actions of a certain central extension $\tJac{N}$ of $\Jac{N}$ by the additive group $\MatT{N}{\RR}$.  It will be necessary for us to work with $\tJac{N}$ in so far as we will use its Casimir element to construct for each slash action an \ido, the {\em Casimir operator.\/}  

\begin{definition}
\label{def:realjacobicentralextension}
Maintaining the $\Mat{N, 2}{\RR}$ identification, the centrally extended rank~$N$ real Jacobi group $\tJac{N}$ and its product are
\begin{gather*}
   \tJac{N} := \bigl\{ (M, X, \kappa):\ (M, X) \in \Jac{N},\
   \kappa \in \Mat{N}{\RR},\ \kappa + \oh X J_2 X^\T \in \MatT{N}{\RR} \bigr\}
\text{,}
\\[6pt]
   (M, X, \kappa) (\check M, \check X, \check \kappa) := 
   (M \check M, X \check M + \check X, 
   \kappa + \check \kappa - X \check M J_2 \check X^\T)
\text{.}
\end{gather*}
\end{definition}

Note that $\Jac{N}$ is centerless, and the center of $\tJac{N}$ is $\MatT{N}{\RR}$.  As we will see in Section~\ref{sec:casimiroperators}, $\tJac{N}$ is a subgroup of $\Sp_{2N+2}{(\RR)}$.

Now fix an element $M := \bigl(\begin{smallmatrix} a & b \\ c & d\end{smallmatrix} \bigr)$ of $\SL{2}(\RR)$.  For $\tau \in \HS$, define
$$ M \tau := (a \tau + b) (c \tau + d)^{-1}, \qquad
   \beta(M, \tau) := (c \tau + d)^{-1}. $$
Then $\tau \mapsto M \tau$ is the standard left action of $\SL{2}(\RR)$ on $\HS$, and $\beta$ is a {\em scalar cocycle\/} with respect to it:
$$ \beta(M \check M, \tau) = \beta(M, \check M \tau) \beta(\check M, \tau). $$
Scalar cocycles are in bijection with slash actions on scalar functions.  For example, $\beta^k$ is a cocycle for all $k \in \ZZ$, and the associated slash action of $\SL{2}(\RR)$ on $C^\infty(\HS)$ is usually written
$$ \phi|_k[M] (\tau) := \beta^k(M, \tau) \phi(M \tau). $$

For future reference, let us mention that the algebra of \dog s on $C^\infty(\HS)$ invariant with respect to the $|_k$-action is the polynomial algebra on one variable generated by the $|_k$-Casimir operator of $\SL{2}(\RR)$, which differs by an additive constant from the weight~$k$ hyperbolic Laplacian
\begin{equation}
\label{eq:laplacian}
   \Delta_k\ :=\ 4y^2 \partial_\tau \partial_{\ov\tau} -
   2iky \partial_{\ov\tau}.
\end{equation}

The theory of cocycles is well-known; see {\em e.g.\/} \cite{BCR10} for a brief summary.  Here we will only review the method by which the scalar cocycles of a given action are classified up to cohomological equivalence.  The stabilizer of $i \in \HS$ under $\SL{2}(\RR)$ is $\SO_2$, and one checks that the restriction of any cocycle to $\SO_2 \times \{ i \}$ defines a representation of $\SO_2$ on $\CC$.  Moreover, it is a fact that two cocycles are equivalent \iff\ they define equal representations of $\SO_2$.  It follows that $\{ \beta^k: k \in \ZZ \}$ exhausts the cocycles of the action under consideration up to equivalence.  For example, the conjugate $\ov\beta\,^k$ is also a cocycle, equivalent to $\beta^{-k}$.

Henceforth write $X_1$ and $X_2$ for the columns of any element $X$ of $\Mat{N, 2}{\RR}$.  The action of $\SL{2}(\RR)$ on $\HS$ generalizes to the following well-known left action of $\Jac{N}$ on $\HS_{1,N}$:
\begin{equation}
\label{eq:definition_HSaction}
  (M, X) (\tau, z) := \bigl(M \tau, \beta(M, \tau) (z + X_1 \tau + X_2) \bigr).
\end{equation}
Regard this as an action of $\tJac{N}$.  As such, the stabilizer of the element $(i, 0)$ of $\HS_{1,N}$ is $\tJacK{N} := \SO_2 \times \{0\} \times \MatT{N}{\RR}$, and the equivalence classes of the scalar cocycles of the action are in bijection with the \rep s of $\tJacK{N}$ on $\CC$.

In order to describe a complete family of cocycles, define a function $a: \tJac{N} \times \HS_{1,N} \to \MatT{N}{\CC}$ by
\begin{eqnarray*}
   a \bigl( (M, X, \kappa), (\tau, z) \bigr) &:=&
   \kappa + X_2 X_1^\T + X_1 z^\T + z X_1^\T + X_1 X_1^\T \tau \\[6pt]
   && - c \beta(M, \tau) (z + X_1 \tau + X_2) (z + X_1 \tau + X_2)^\T
\end{eqnarray*}
(recall that $c$ is $M_{21}$).  For $L \in \MatT{N}{\CC}$, define $\alpha_L: \tJac{N} \times \HS_{1,N} \to \CC$ by
$$ \alpha_L \bigl( (M, X, \kappa), (\tau, z) \bigr)\ :=\
   \exp \bigl\{ 2 \pi i \thinspace \tr 
   \bigl[L a \bigl( (M, X, \kappa), (\tau, z) \bigr) \bigr]\bigr\}. $$

\begin{lemma} \label{lemma:cocycles}
For all $k \in \ZZ$ and $L \in \MatT{N}{\CC}$, $\beta^k \alpha_L$ is a scalar cocycle with respect to the action~(\ref{eq:definition_HSaction}) on $\HS_{1,N}$ of the centrally extended Jacobi group $\tJac{N}$ from Definition~\ref{def:realjacobicentralextension}.  Moreover, any scalar cocycle of this action is equivalent to exactly one of these cocycles.
\end{lemma}

\begin{proof}
The proof that $\beta^k$ is a cocycle of the action of $\tJac{N}$ on $\HS_{1,N}$ is the same as the proof that it is a cocycle of the action of $\SL{2}(\RR)$ on $\HS$.  The proof that $\alpha_L$ is a cocycle is standard in the case $N = 1$ and proceeds along the same lines in general.  One must prove that $a(g \check g, x) = a(g, \check g x) + a(\check g, x)$.  First check that it suffices to prove this for both $g$ and $\check g$ in either the semisimple or the nilpotent part of $\tJac{N}$, and then check each of the resulting four cases directly.  The second sentence follows immediately from the classification of \rep s of $\tJacK{N}$.
\end{proof}

As a consequence of this lemma we have the following family of slash actions of $\tJac{N}$ on $C^\infty(\HS_{1,N})$: for $k$, $k' \in \ZZ$ and $L \in \MatT{N}{\CC}$, 
\begin{eqnarray*}
  \phi|_{k, k', L} [M, X, \kappa] (\tau, z) &:=&
  \phi \bigl((M, X, \kappa)(\tau,z)\bigr) \\[6pt]
  && \times\thinspace \beta^k (M, \tau) {\overline\beta}\,\!^{k'} (M, \tau)
  \alpha_L \bigl((M, X, \kappa)(\tau,z)\bigr).
\end{eqnarray*}
Observe that since $\beta \ov\beta$ is positive, $|_{k,k',L}$ makes sense for all $k$, $k' \in \CC$ with $k - k' \in \ZZ$.  We will write $|_{k,L}$ for $|_{k,0,L}$.  (Usually we will be concerned only with the case $k' = 0$, but at one point we will need the freedom to choose differently.)  By Lemma~\ref{lemma:cocycles}, any slash action is equivalent to exactly one of the actions $|_{k,L}$; as we have mentioned, $|_{k,k',L}$ is equivalent to $|_{k-k',L}$.  

\begin{definition} \label{defn:CDOs}
A \dog\ $T$ on\/ $\HS_{1,N}$ is\/ {\em covariant from $|_{k,L}$ to $|_{k',L'}$\/} if for all $g \in \tJac{N}$ and $f \in C^\infty(\HS_{1,N})$, we have
$$ T \bigl( f|_{k,L}[g] \bigr)\ =\ (T f) \big|_{k',L'}[g]. $$
Let $\bD(k,L;\, k',L')$ be the space of covariant operators from $|_{k,L}$ to $|_{k',L'}$, and let $\bD^r(k,L;\, k',L')$ be the space of those of order~$\le r$.  When $k' = k$ and $L' = L$, we refer to such operators as\/ {\em $|_{k,L}$-invariant\/} and write simply\/ $\bD_{k,L}$ and\/ $\bD^r_{k,L}$.
\end{definition}

At this point we state the main results of Section~\ref{sec:casimiroperators}, Theorem~\ref{thm:actionofcenter} and Propositions~\ref{prop:RLops}, \ref{prop:allCDOs}, and~\ref{prop:allIDOs}.  Elements of $C^\infty(\HS_{1,N})$ holomorphic in $\bC^N$ will be called {\em semi-holomorphic.\/}  For any $N \times N$ matrix $A$ and any $N$-vector $w$, set
$$ A[w] := w^\T A w .$$
Recall the Laplacian~(\ref{eq:laplacian}) and our notation $\tau := x + iy \in \bC$ and $z := u + iv \in \bC^N$.  For brevity, write $\sL := 2\pi i L$.  For $L$ invertible, define
\begin{eqnarray}
   \cC^{k,L} & := &
   - 2 \Delta_{k-N/2} 
   + 2 y^2 \bigl( \partial_{\ov\tau} \sL^{-1}[\pz] + \partial_\tau \sL^{-1}[\pbz] \bigr)
   - 8y \partial_\tau v^\T \pbz \nonumber\\[6pt]
&& - {\ts\frac{1}{2}} y^2 \bigl( \sL^{-1}[\pbz] \sL^{-1}[\pz] - (\pbz^\T \sL^{-1} \pz)^2 \bigr)
   + 2 y (v^\T \pbz) \pz^\T \sL^{-1} \partial_u \label{casop}\\[6pt]
&& - {\ts\frac{1}{2}} (2k - N + 1) iy \pbz^\T \sL^{-1} \partial_u
   + 2 v^\T (v^\T \pbz) \pbz
   + (2k - N -1) iv^\T \pbz. \nonumber
\end{eqnarray}

\begin{theorem}
\label{thm:actionofcenter}
For $L$ invertible, the operator\/ $\cC^{k,L}$ is, up to additive and multiplicative scalars, the {\em Casimir operator\/} of\/ $\tJac{N}$ with respect to the\/ $|_{k, L}$-action (see Section~\ref{sec:casimiroperators}).  It generates the image of the\/ $|_{k,L}$-action of the center of the \uea\ of\/ $\tJac{N}$.  In particular, it lies in the center of\/ $\bD_{k,L}$.  Its action on semi-holomorphic functions is
\begin{align}
\label{eq:semiholocasimiroperator}
  -2 \Delta_{k - N/2} + 2 y^2 \partial_{\ov\tau} \sL^{-1}[\pz].
\end{align}
\end{theorem}

Note that for $N>1$, (\ref{casop})~is of order~4.  At $N=1$ it is of order~3 and reduces to the operator $C^{k,m}$ given in \cite{BR10} with $L=m$.  (There is a misprint in \cite{BR10}: the term $k (z - \ov z) \pbz$ should be $(1 - k) (z - \ov z) \pbz$.  This stems in part from a similar misprint in~(8) of \cite{Pi09}, where the term $(z - \ov z) \pbz$ coming from (6) of \cite{Pi09} is missing.)

\begin{definition} \label{defn:RLops}
The\/ {\em lowering operators,\/} $X_-$ and\/ $Y_-$, and the\/ {\em raising operators,\/} $X_+$ and\/ $Y_+$, are
\begin{align*}
   X^{k,L}_- &:= -2iy \bigl(y \partial_{\ov\tau} + v^\T \pbz \bigr),\ &
   X^{k,L}_+ &:= 2i \bigl(\partial_\tau + y^{-1} v^\T \pz + y^{-2} \sL[v]\bigr)
   + k y^{-1}, \\[6pt]
   Y^{k,L}_- &:= -iy \pbz,\ &
   Y^{k,L}_+ &:= i \pz + 2i y^{-1} \sL v.
\end{align*}
\end{definition}

For $N=1$ and $L=m$, these are the operators given on page~59 of \cite{BS98}.  (There is a misprint in their formula for $Y_-$: the expression $\oh (\tau - \ov\tau) f_{\ov z}$ on the far right should be multiplied by $-1$.)  Note that $Y^{k,L}_\pm$ are actually $N$-vector operators.  We write $Y^{k,L}_{\pm, j}$ for their entries.

Frequently we will suppress the superscript $(k,L)$.  Care must be taken with this abbreviation, as for example $X_+ Y_+$ means $X^{k+1, L}_+ Y^{k,L}_+$.

\begin{proposition} \label{prop:RLops}
The spaces\/ $\bD^1(k,L;\, k \pm 2, L)$ are 1-dimensional, and the spaces\/ $\bD^1(k,L;\, k \pm 1, L)$ are $N$-dimensional.  They have bases given by
\begin{equation*}
   \bD^1(k,L;\, k \pm 2, L) = \Span \bigl\{ X^{k,L}_\pm \bigr\}, \quad
   \bD^1(k,L;\, k \pm 1, L) = \Span \bigl\{ Y^{k,L}_{\pm, j}:\, 1 \le j \le N \bigr\}.
\end{equation*}
The spaces $\bD^1_{k, L}$ are equal to $\bD^0_{k, L} = \bC$.  All other\/ $\bD^1(k,L;\, k',L')$ are zero.

The raising operators commute with one another, as do the lowering operators (but keep in mind that, for example, $X_+ Y_+ = Y_+ X_+$ means $X_+^{k+1, L} Y_+^{k, L} = Y_+^{k+2, L} X_+^{k, L}$).  The commutators between them are
\begin{equation*}
   [X_-,\, X_+] = -k, \quad
   [Y_{-,j},\, Y_{+,j'}] = i \sL_{jj'}, \quad
   [X_-,\, Y_+] = -Y_-, \quad
   [Y_-,\, X_+] = Y_+.
\end{equation*}
\end{proposition}

\begin{proposition} \label{prop:allCDOs}
Any covariant \dog\ of order~$r$ may be expressed as a linear combination of products of up to $r$ raising and lowering operators.  There is a unique such expression in which the raising operators are all to the left of the lowering operators.

The expression of this form for the Casimir operator is
\begin{equation}
\begin{array}{rcl}
\label{eq:casimir_raisinglowering}
   \cC^{k,L} & = &
   -2 X_+ X_- + i\bigl( X_+ \sL^{-1}[Y_-] - \sL^{-1}[Y_+] X_- \bigr) \\[6pt]
&& -\oh \bigl( \sL^{-1}[Y_+] \sL^{-1}[Y_-]
   - Y_+^\T (Y_+^\T \sL^{-1} Y_-) \sL^{-1} Y_- \bigr) \\[6pt]
&& - \oh (2k - N - 3) i Y_+^\T \sL^{-1} Y_-.
\end{array}
\end{equation}
\end{proposition}

\begin{proposition} \label{prop:allIDOs}
The algebra\/ $\bD_{k,L}$ is generated by\/ $\bD^3_{k,L}$.  The spaces\/ $\bD^3_{k,L}$ and\/ $\bD^2_{k,L}$ are of dimensions $2N^2 + N + 2$ and $N^2 + 2$, respectively.  Bases for them are given by the following equations:
\begin{eqnarray*}
   \bD^3_{k,L} &=&
   \Span \bigl\{ X_+ Y_{-,i} Y_{-,j},\ Y_{+,i} Y_{+,j} X_-:\, 1 \le i \le j \le N \bigr\}
   \, \oplus\, \bD^2_{k,L}, \\[6pt]
   \bD^2_{k,L} &=&
   \Span \bigl\{ 1,\ X_+ X_-,\ Y_{+, i} Y_{-, j}:\, 1 \le i, j \le N \bigr\}.
\end{eqnarray*}
\end{proposition}

The focus of this paper is the space of harmonic Maa\ss-Jacobi forms of index $L$ and weight $k$.  In order to define it, fix $k \in \ZZ$ and a positive definite integral even lattice $L$ of rank $N$.  We will identify $L$ with its Gram matrix with respect to a fixed basis, a positive definite symmetric matrix with entries in $\frac{1}{2} \ZZ$ and diagonal entries in $\ZZ$.  Write $|L|$ for the covolume of the lattice, the determinant of the Gram matrix.

The full Jacobi group $\JacF{N}$ defined in~(\ref{eq:jacandfulljac}) clearly has a central extension by $\MatT{N}{\ZZ}$ which is a subgroup of $\tJac{N}$.  It is easy to check that when $L$ is a Gram matrix, the cocycle $\alpha_L$ is trivial on $\MatT{N}{\ZZ}$.  Therefore the $|_{k,L}$-action factors through to an action of $\JacF{N}$, which we will also denote by $|_{k,L}$.

\begin{definition}[Maa\ss-Jacobi forms]
\label{def:maassjacobiforms}
A Maa\ss-Jacobi form of weight $k$ and index $L$ is a function $\phi \in C^{\infty} (\HS_{1,N})$ satisfying the following conditions:
\begin{enumerate}
\item For all $A \in \JacF{N}$, we have $\phi|_{k, L} [A] = \phi$.
\item $\phi$ is an eigenfunction of\/ $\cC^{k,L}$.
\item For some $a > 0$, $\phi(\tau, z) = O\left( e^{a y} e^{2 \pi L[v] / y} \right)$ as $y \rightarrow \infty$.
\end{enumerate}
If $\phi$ is annihilated by the Casimir operator $\cC^{k,L}$, it is said to be a harmonic Maa\ss-Jacobi form.  We denote the space of all harmonic Maa\ss-Jacobi forms of fixed weight $k$ and index $L$ by $\bJ_{k,L}$.
\end{definition}

\begin{remark}
Adapting the proof in \cite[Section 2.6]{BS98}, which is based on \cite[Section 1.3]{LV80} and \cite[Section 2.I.2]{MVW87}, we see that any automorphic representation of $\tJac{N}$ is a tensor product $\tilde \pi \otimes \pi^L_{\mathrm{SW}}$.  Here $\tilde \pi$ is a genuine representation of the metaplectic cover of\/ $\SL{2}$, and $\pi^L_{\mathrm{SW}}$ is the Schr\"odinger-Weil representation of central character $L$.  The latter is the extension to the metaplectic cover of the Jacobi group of the Schr\"odinger representation of the Heisenberg group, which is induced from the character $e^{2 \pi i \tr (L \kappa)}$ of its center.  Thus, as in \cite{Pi09}, semi-holomorphic forms play an important role in the representation-theoretic treatment of harmonic Maa\ss-Jacobi forms.
\end{remark}

For later use we set 
$$ e(r) := e^{2 \pi i r}, \
   q := e(\tau), \
   \zeta^r := \prod_{i = 1}^N e(z_i r_i)
\text{.}
$$


%% file: mockjacobiforms.tex
\section{H-harmonic Jacobi forms}
\label{sec:mockjacobiforms}

The concept of H-harmonic (Heisenberg harmonic) Jacobi forms for $\JacF{1}$ was introduced in \cite{BRRZ11}: they define a function $\phi \in C^\infty(\HS_{1,1})$ to be H-harmonic if $Y_+ Y_- \phi = 0$.  Imposing H-harmonicity on harmonic Maa\ss-Jacobi forms amounts to fixing a 4-dimensional space of Fourier coefficients for every index.  This is what makes this notion as important as it is.  H-harmonicity incorporates famous examples of real-analytic Jacobi forms, such as the Zwegers $\mu$-function.  It also allows for a decomposition similar to the classical $\theta$-decomposition.

In this section we will generalize the notion of H-harmonicity to all Jacobi groups $\tJac{N}$.  The main result is that for $N \ge 2$, H-harmonicity coincides with the notion of semi-holomorphicity investigated in Section~\ref{sec:semiholomorphicforms}.

There are three fundamental requirements which any definition of H-harmonicity in $C^\infty(\HS_{1,N})$ should satisfy.  First, it should be independent of coordinates, because \`a priori there are no distinguished directions in the second component of $\HS_{1,N}$.  Second, specializing an H-harmonic function to $\HS_{1,1}$ by restricting $z$ to the line through a fixed $v \in \QQ^N$ should give an H-harmonic function.  Third, the definition should be based on covariant operators, because we will ultimately consider H-harmonic Jacobi forms which are invariant under $\JacF{N}$.  Imposing vanishing with respect to differential operators that are not covariant automatically leads to a stronger vanishing condition with respect to covariant operators.

Specializing the second component of an H-harmonic function $\phi$ to $(z_1, 0, \ldots)$, the second and third requirements lead to the natural condition $Y_{+,1} Y_{-,1} \phi = 0$.  Symmetrizing this to meet the first requirement we obtain the following definition.
\begin{definition}
\label{def:hharmonic}
A function $\phi \in C^{\infty} (\HS_{1,N})$ is said to be {\em H-harmonic} if it vanishes under $(B Y_+)^\T (B Y_-)$ for all $B \in \GL{N}(\RR)$.
\end{definition}

We connect this definition with invariant metrics on $\HS_{1,N}$.  To make their expression more readable we use the {\em S-coordinates\/} $(p, q)$ on $\bC^n$ defined by $z = p\tau + q$ with $p$ and $q$ real (see \cite{BS98}).

\begin{proposition}
\label{prop:laplaceoperators}
For any positive definite symmetric matrix\/ $C \in \MatT{N}{\RR}$,
\begin{align*}
  ds^2
=
   y^{-2} \partial_\tau \partial_{\overline{\tau}}
   + y^{-1} (\tau \overline{\tau}\, C[\partial_p] + 2 x\, \partial_q^\T C \partial_p + C[\partial_q])
\end{align*}
is an invariant metric.  The associated Laplace operator is $X_+ X_- + Y_+^\T C Y_-.$
\end{proposition}

\begin{proof}
The invariance with respect to $\SL{2}(\RR) \subseteq \Jac{N}$ follows as for $N = 1$.  The invariance with respect to the Heisenberg group can be proven by choosing an appropriate basis of $\CC^N$ and again following the calculation for $N = 1$.  The associated Laplace operator can be computed by choosing a basis of $\CC^N$ with respect to which $C$ is diagonal.
\end{proof}

A function $\phi \in C^{\infty} (\HS_{1,N})$ is harmonic with respect to all Laplace operators in Proposition \ref{prop:laplaceoperators} if and only if it vanishes under the operators $X_+ X_-$ and $Y_{+} C Y_{-}$ for all (not necessarily invertible) $C \in \MatT{N}{\RR}$.  Note that the second condition is equivalent to the condition imposed in Definition~\ref{def:hharmonic}.

\begin{proposition}
\label{prop:restriction_of_hharmonic}
Fix $\phi \in C^\infty (\HS_{1,N})$ as above.  For any $N' < N$ and any matrix $M \in \Mat{N,N'}{\QQ}$ with full rank, the function $(\tau, z') \mapsto \phi(\tau, M z')$ from $\HS_{1,N'}$ to $\CC$  is H-harmonic whenever $\phi$ is.
\end{proposition}
\begin{proof}
The statement reduces to linearity of differential operators.
\end{proof}

Using this proposition, we now prove that H-harmonic Jacobi forms are not more general than semi-holomorphic forms for $N \ge 2$.  Hence we may say that the Heisenberg part of $\tJac{N}$ is more rigid in the higher rank case.

\begin{theorem}
If $\phi \in C^\infty(\HS_{1,N})$ is an H-harmonic Maa\ss-Jacobi form and $N \ge 2$, then $\phi$ is holomorphic in $z$.
\end{theorem}
\begin{proof}
It suffices to proof that $\phi$ is holomorphic in every coordinate of $z$.  Thus, by Proposition~\ref{prop:restriction_of_hharmonic}, we may restrict to the case $N = 2$.  Since H-harmonicity is coordinate-free, we may fix coordinates so that $L_{12} \ne 0$.  With these coordinates fixed, we prove that if $Y_{+,i} Y_{-,i} \phi = 0$ for $i = 1,2$, then $\phi$ is semi-holomorphic.  For reasons of symmetry, it suffices to show that $\partial_{\ov z_1}\, \phi = 0$, or equivalently, $Y_{-,1} \phi = 0$.

Using the Fourier expansion of $\phi$ it is sufficient to consider a single term $a(y, v) q^n \zeta^r$.  Since the $Y_{\pm}$ are covariant, we may assume that $r = 0$.  We are reduced to showing that any smooth function annihilated by
\begin{equation*}
  A = \partial_{v_1}^2 + 2 (L_{11} v_1 + L_{12} v_2) \partial_{v_1}
\quad \text{and} \quad
  B = \partial_{v_2}^2 + 2 (L_{12} v_1 + L_{22} v_2) \partial_{v_2}
\end{equation*}
is already a constant.  Direct verification shows that
\begin{equation*}
  \partial_2 - \partial_1
=
  \frac{1}{2 L_{12}}
  \Big( \partial_{v_2}^2 A - \partial_{v_1}^2 B
        - (L_{12} v_1 + L_{22} v_2) \partial_{v_2} A
        + (L_{11} v_1 + L_{12} v_2) \partial_{v_1} B
  \Big)
\text{.}
\end{equation*}
Using this relation, we deduce from the vanishing under $A$ and $B$ that $a(v)$ is annihilated by
\begin{eqnarray*}
  \partial_{v_1}^2 + 2 (L_{11} v_1 + L_{12} v_2) \partial_{v_1}
  -
  \partial_{v_1}^2 - 2 (L_{12} v_1 + L_{22} v_2) \partial_{v_1} &&
\\[6pt] \ =\
  2 ((L_{11} - L_{12}) v_1 - (L_{22} - L_{12}) v_2) \partial_{v_1}
\text{.} &&
\end{eqnarray*}
In other words, $a(v)$ is constant outside of a 1-dimensional space.  Since it is smooth, it is constant everywhere.
\end{proof}


%% file: semiholomorphicforms.tex
\section{Semi-holomorphic forms}
\label{sec:semiholomorphicforms}

Recall that a function on $\HS_{1,N}$ holomorphic in $z \in \CC^N \subseteq \HS_{1,N}$ is called semi-holomorphic.  We will denote the space of semi-holomorphic harmonic Maa\ss-Jacobi forms by $\mathbb{J}_{k,L}^{\mathrm{semi}}$.  Semi-holomorphic forms vanish under $Y_-$, and $X_-$ acts on them by $\partial_{\overline{\tau}}$.  In particular, they fall under Definition \ref{def:hharmonic}.

The theory of semi-holomorphic forms essentially mimics that of harmonic weak Maa\ss\ Forms.  Indeed, in Theorem \ref{thm:thetadecompositioncompatibility} we will see that the $\theta$-decomposition gives a well-behaved bijection between vector-valued weak harmonic Maa\ss\ forms and harmonic semi-holomorphic Maa\ss-Jacobi forms.

We first discuss semi-holomorphic Fourier expansions of Maa\ss-Jacobi forms.  The negative discriminant of a Fourier index $(n, r)$ is denoted by
\begin{gather*}
  D
:=
  D_L(n,r) := |L| ( 4 n - L^{-1}[r])
 \end{gather*}
By analogy with \cite[page 9]{BF04}, define a function
\begin{align*}
  H(y)
& :=
  e^{-y} \int_{-2 y}^{\infty} e^{-t} t^{-k - N / 2} \;dt
\text{.}
\end{align*}

\begin{proposition}
\label{prop:fourierexpansionofmaassjacobiforms}
Any semi-holomorphic harmonic Maa\ss-Jacobi form $\phi$ has a Fourier expansion of the form
\begin{align*}
  \phi(\tau, z)
=
  y^{N/2 - k} &\!\!\sum_{\substack{n \in \ZZ, r \in \ZZ^N \\ \text{s.t. } D = 0}}
              \!\!c^0(n,r) \; q^n \zeta^r
+ \sum_{\substack{n \in \ZZ, r \in \ZZ^N \\ \text{s.t. } D \gg -\infty}}
  \!\!c^+ (n,r) \; q^n \zeta^r
\\&
  + \sum_{\substack{n \in \ZZ, r \in \ZZ^N \\ \text{s.t. } D \ll \infty}}
  \!\!c^- (n,r) \, H(\pi D y \slashdiv 2 |L|) \, e\bigl(-i D y \slashdiv 4 |L|\bigr)\; q^n \zeta^r,
\end{align*}
where the Fourier coefficients $c^0(n,r)$ and $c^\pm(n,r)$ are in $\CC$.
\end{proposition}

\begin{proof}
This can be proved as in the case of rank $1$ lattices, by solving the differential equation for the coefficients coming from the Casimir operator and then imposing the growth condition.
\end{proof}

Our investigation will concentrate on semi-holomorphic harmonic Maa\ss-Jacobi forms, and in particular their relation to skew-holomorphic forms.  To state this relation we must define a $\xi$-operator.  Proceeding as in \cite[Section 4]{BR10}, we first define the lowering operator
\begin{align*}
\label{eq:definition_loweroperator}
  D_{-}^{(L)}
& :=
  -2 i y \bigl( y \, \partial_{\overline \tau}
               + v^\T \partial_{\overline z}
               - {\ts \frac{1}{4}} y \sL^{-1} [\partial_{\overline z}] \bigr)
   \ =\ X_- - {\ts \frac{i}{2}} \sL^{-1}[Y_-]
\text{.}
\end{align*}
Using this operator, we define the $\xi$-operator by
\begin{align}
\label{eq:definition_xioperator}
  \xi_{k,L}
& :=
   y^{k - 2 - N/2} D_{-}^{(L)}
\text{.}
\end{align}
This is an analog of the $\xi$-operator in \cite{Ma49}.  The latter sends Maa\ss\ forms to their shadows, which are holomorphic if they have harmonic preimages.  In our setting skew-holomorphic forms take the place of holomorphic ones.

\begin{definition}[Skew-holomorphic Jacobi forms]
\label{def:skewholomorphicforms}
A skew-holomorphic Jacobi form of weight $k$ and index $L$ is a semi-holomorphic function $\phi \in C^{\infty} (\HS_{1,N})$ satisfying the following conditions.  First, for all $A \in \JacF{N}$ the equation $ \left.\phi \right|_{N/2,k-N/2,L} A = \phi$ holds.  Second, the Fourier expansion of $\phi$ has the form
\begin{align*}
  \phi(\tau, z)
& =
  \sum_{\substack{n \in Z, r \in \ZZ^N \\ \text{s.t. } D \gg -\infty}}
  c(n, r) \, e(- i D y \slashdiv 2 |L|) \; q^n \zeta^r
\text{.}
\end{align*}
We write $\mathbb{J}_{k,L}^{\mathrm{sk}}$ for the space of all such forms.
\end{definition}

\begin{remark}
Skew-holomorphic Jacobi forms were first introduced by Skoruppa in \cite{Sk90}.  There are several articles treating a slightly more general notion than that we have given.  See in particular \cite{Ha06}.
\end{remark}

\begin{remark}
The Fourier expansion condition can be stated in terms of annihilation by the heat operator $2 \partial_{\tau} - \sL^{-1}[\partial_{z}] \slashdiv 2$.
\end{remark}

\begin{proposition}
If $\phi \in \mathbb{J}_{k,L}^{\mathrm{semi}}$, then $\xi_{k,L}\phi$ is an element of $\mathbb{J}_{2 + N - k,L}^{\mathrm{sk}}$.
\end{proposition}

\begin{proof}
By Proposition~\ref{prop:RLops}, $D_{-}^{(L)}$ is a covariant operator from $|_{k,L}$ to $|_{k-2,L}$.  Applying $\xi_{k,L}$ to the Fourier expansion of a Maa\ss-Jacobi form as in Proposition~\ref{prop:fourierexpansionofmaassjacobiforms} shows that the Fourier expansion of $\xi_{k,L} \phi$ has the correct form.
\end{proof}

The $\xi$-operator is compatible with the $\theta$-decomposition.  To state this precisely, let $\Gamma$ be the elliptic metaplectic group with the same level as $L$.  Denote the spaces of vector-valued harmonic Maa\ss\ forms for the Weil representation $\rho_L$ by $[\Gamma, k - N \slashdiv 2, \rho_L]^{\mbox{\scrm Maa\ss}}$.  For weakly holomorphic vector-valued Maa\ss\ forms change the superscript to ${*}^{\mathrm{!}}$.  The $\xi$-operator $\xi_{k - N / 2} = y^{k - N / 2} \, {\overline \partial}_{\overline{\tau} \; \cdot }$ maps this space of harmonic Maa\ss\ forms to the space of weakly holomorphic forms.

To revise the $\theta$-decomposition we need the following $\theta$-series for $\mu \in \ZZ^N$:
\begin{align}
\label{eq:thetaseriesL_definition}
  \theta_{L,\mu} (\tau, z)
& :=
  \sum_{r \in \ZZ^N,\, r \equiv \mu (L \ZZ^N)} \!\! q^{L^{-1}[r] / 4} \zeta^r
\text{.}
\end{align}
 
\begin{definition}[$\theta$-decomposition]
The Maa\ss-Jacobi $\theta$-decomposition is the map $\theta_L^{\mathrm{semi}} : \mathbb{J}_{k,L}^{\mathrm{semi}} \rightarrow [\Gamma, k - N \slashdiv 2, \rho_L]$ defined by
\begin{align*}
  f(\tau, z)
& =
  \sum_{\mu (\ZZ^N / L \ZZ^N)} \theta_{L}^{\mathrm{semi}}(f)_\mu (\tau) \; \theta_{L, \mu} (\tau, z)
\text{.}
\end{align*}
Similarly, the skew-holomorphic $\theta$-decomposition map $\theta_{L}^{\mathrm{sk}} : \mathbb{J}_{k,L}^{\mathrm{sk}} \rightarrow [\Gamma, k - N \slashdiv 2, \rho_L]$ is defined by
\begin{align*}
  f(\tau, z)
& =
  \sum_{\mu (\ZZ^N / L \ZZ^N)} {\overline \theta}_{L}^{\mathrm{sk}}(f)_\mu (\tau) \; \theta_{L, \mu} (\tau, z)
\text{.}
\end{align*}
\end{definition}
\begin{remark}
The existence of a $\theta$-decomposition for a harmonic Maa\ss\ form is equivalent to its semi-holomorphicity.
\end{remark}

\begin{theorem}
\label{thm:thetadecompositioncompatibility}
If $k$ is even, the $\theta$-decomposition of forms in $\mathbb{J}_{k,L}^{\mathrm{semi}}$ and $\mathbb{J}_{2 + N - k,L}^{\mathrm{sk}}$ commutes with the $\xi$-operators $\xi_{k,L}$ and $\xi_{k - N \slashdiv 2}$.  More precisely, the following diagram is commutative:
\begin{align*}
\xymatrixcolsep{4pc}
\xymatrix{
 \mathbb{J}_{k,L}^{\mathrm{semi}} \ar[d]^{\theta_{L}^{\mathrm{semi}}} \ar[r]^{\xi_{k,L}} &
 \mathbb{J}_{2 + N - k,L}^{\mathrm{sk} }\ar[d]^{\theta_L^{\mathrm{sk}}}\\
 \lbrack \Gamma, k - N \slashdiv 2, \rho_L\rbrack^{\mbox{\scrm Maa\ss}} \ar[r]^-{\xi_{k - N \slashdiv 2}} &
 \lbrack\Gamma, 2 + N \slashdiv 2 - k, \rho_L\rbrack^{\mathrm{!}} 
}
\text{.}
\end{align*}
\end{theorem}

\begin{proof}
This is a calculation analogous to that in \cite[Section 6]{BR10}.
\end{proof}

Before we consider the Poincar\'e series we define a special part of the space of semi-holomorphic harmonic Maa\ss-Jacobi forms. We will show that it maps surjectively to the space of skew-holomorphic Jacobi forms with cuspidal shadow.

\begin{definition}[Maa\ss-Jacobi forms with cuspidal shadow]
The inverse image under $\xi_{k,L}$ of\/ $\mathbb{J}_{k,L}^{\mathrm{sk},\mathrm{cusp}}$, the cuspidal subspace of $\mathbb{J}_{k,L}^{\mathrm{sk}}$, is denoted by $\mathbb{J}_{k,L}^{\mathrm{semi},\mathrm{cusp}}$.  It is the space of semi-holomorphic harmonic Maa\ss-Jacobi forms with cuspidal shadow.
\end{definition}

\subsection{Poincar\'e series}
\label{ssec:poincareseries}

In \cite[Section 5]{BR10} the authors define Maa\ss-Poincar\'e series for the Jacobi group.  They restrict to Jacobi indices of rank one.  In this section we generalize their considerations to arbitrary lattice indices.

We use the notation of Section~\ref{sec:maassjacobiforms}; in particular, $L$ is an integral lattice and $k$ is in $\ZZ$.  Throughout this section $n$ will be an integer and $r$ will be in $\ZZ^N$.  Maintain $D$ as above and set $h$ as follows:
$$ D := D_L(n, r) := |L|(4 n - L^{-1}[r]), \qquad h := h_L(r) := |L|\, L^{-1}[r]. $$
The standard scalar product of two $N$-vectors $\lambda$ and $z$ will be written as $\lambda z$.

Using the $M$-Whittaker function $M_{\nu,\mu}$ (see \cite{WW96}), we define
\begin{align}
\label{eq:whittakerM_definition}
  \wtM_{s,\kappa} (t)
& :=
 |t|^{- \kappa / 2} M_{\mathrm{sgn}(t) \kappa / 2,\, s - 1/2}(|t|)
\text{,}\\[6pt]
\label{eq:phikls_definition}
  \phi_{k, L, s} ^{(n, r)} (\tau, z)
& :=
  \wtM_{s,k - N/2} (\pi D y \slashdiv |L|) \;
  e(r z + i L^{-1}[r] y \slashdiv 4 + n x)
.
\end{align}

\begin{lemma}
The function $\phi_{k,L,s}^{(n,r)}$ defined in~(\ref{eq:phikls_definition}) is an eigenfunction of the Casimir operator $\cC^{k,L}$ in Theorem~\ref{thm:actionofcenter}, with eigenvalue
\begin{gather} \label{eq:eval}
  -2 s ( 1 - s ) - {\textstyle\frac{1}{2}} \Bigl(k^2 - k (N + 2) + {\textstyle\frac{1}{4}} N (N + 4) \Bigr)
\text{.}
\end{gather}
\end{lemma}

\begin{proof}
Factor $\phi$ as follows:
\begin{align*}
  \phi_{k,L,s}
& =
    e(r z + \tau L^{-1}[r] \slashdiv 4)
\;\cdot\;
    e(-D x \slashdiv 4 |L|) \, \wtM_{s,k - N/2} (- \pi D y \slashdiv |L|)
\text{.}
\end{align*}
The first factor is holomorphic in $\tau$ and the second is constant in $z$.  Hence in applying $\cC^{k,L}$ the contribution of the first factor cancels.  We need only consider $-2 \Delta_{k - N/2}$, yielding~(\ref{eq:eval}).
\end{proof}

We will study the Poincar\'e series
\begin{align}
\label{eq:poincareseries_definition}
  P_{k,L,s}^{(n,r)}
:=
  \sum_{A \in \JacFpara{N} \backslash \JacF{N}} \left.\phi_{k,L,s}^{(n,r)} \right|_{k,L} A
,
\end{align}
which is semi-holomorphic.  The usual estimate
\begin{align*}
  \wtM_{s,k - N/2} (y)
\ll
  y^{\Re(s) - (2 k - N)/4}
\quad\text{as}\quad
  y \rightarrow 0
\end{align*}
ensures absolute and uniform convergence for $\Re(s) > 1 + N \slashdiv 2$.  Of particular interest will be the case $s \in \{k/2 - N/4,\,1 + N/4 - k/2\}$, where the Poincar\'e series is annihilated by the Casimir operator.
\vspace*{1em}

We need to compute the Fourier expansions of the Poincar\'e series, which involve the $I$-Bessel function as well as the $J$-Bessel function.  The following $W$-Whittaker function, $\theta$-series, and higher Kloosterman sum will also arise.  To make the notation more natural we renormalize the Whittaker function:
\begin{align}
\label{eq:whittakerW_definition}
  \wtW_{s,\kappa} (t)
& :=
  |t|^{-\kappa / 2} \;W_{\mathrm{sgn}(t) \kappa / 2,\, s - 1/2}(|t|)
\text{,}\\[6pt]
\label{eq:thetaseries_definition}
  \theta_{k,L}^{(r)}
& := 
  \sum_{\lambda \in \ZZ^N} q^{L[\lambda]} \zeta^{2 L \lambda} \,
  (q^{r \lambda} \zeta^r + (-1)^k q^{-r \lambda} \zeta^r)
\text{,}\\[6pt]
\label{eq:klostermansum_definition}
  K_{c,L}(n,r,n',r')
& :=
  e(-r L^{-1} r' \slashdiv 2 c)
\\[6pt] \nonumber
  & \hphantom{:=}
  \sum_{d (c)^\times,\,\lambda \in \ZZ^N \slashdiv c \ZZ^N} \hspace{-0.3em}
  e(\overline{d} L[\lambda] \slashdiv c + n' d - r' \lambda
    + \overline{d} n + \overline{d} r \lambda)
\text{,}
\end{align}
where $\overline{d}$ is an integer inverse of $d$ modulo $c$.

\begin{theorem}
\label{th:poincarefourierexpansions}
The Poincar\'e series~(\ref{eq:poincareseries_definition}) has the Fourier expansion
\begin{align}
\label{eq:poincarefourierexpansion}
  P_{k, L, s}^{(n,r)} (\tau, z)
& =
  \wtM_{s, k - N/2}(\pi D y \slashdiv |L|) \, e(- i D y \slashdiv 4 |L|)\;
  \theta_{k, L}^{(r)} (\tau, t) q^n
\\&
\nonumber
  + \sum_{n' \in \ZZ, r' \in \ZZ^N} \!c_{y,s}(n', r') \; q^{n'} \zeta^{r'}
\text{.}
\end{align}
Here the $\theta$-series $\theta_{k, L}^{(r)}$ is defined in~(\ref{eq:thetaseries_definition}), and the coefficients $c_{y,s}$ are
\begin{align*}
  c_{y,s}(n', r')
& :=
  b_{y,s}(n', r') + (-1)^k b_{y,s}(n', -r'),
\end{align*}
with $b_{y,s}$ depending on $D$ and $D' = |L|(4 n' + L^{-1}[r'])$.

For $D' = 0$, there is a constant $a_s(n', r')$ such that
\begin{align*}
  b_{y,s} (n', r')
& =
  \frac{y^{1 + N/4 - k/2 - s}}{\Gamma(s + k/2 - N/4)\, \Gamma(s - k/2 + N/4)} \;
  a_s(n', r').
\end{align*}
For $D' \ne 0$,
\begin{align*}
  b_{y,s} (n', r')
& =
  2^{1 - N / 2} \pi i^{-k} \, |L|^{-1/2} \; \Bigl(\Gamma(2 s) \slashdiv \Gamma\bigl(s - \mathrm{sgn}(D') (k/2 - N/4)\bigr)\Bigr)
\\[4pt]&
  \cdot (D'/D)^{k/2 - (N+2)/4} \;\, e(-i D' y \slashdiv 4 |L|) \;
  \mathcal{W}_{s, k - N/2}(\pi D' y \slashdiv |L|)
\\[4pt]&
  \cdot \sum_{c \in \ZZ_{>0}} c^{-(N + 2)/2} K_{c,L}(n, r, n', r')
\\&
  \hphantom{\sum_{c \in \ZZ_{>0}} c^{-(N + 2)/2}}
  \cdot
  \begin{cases}
  J_{2 s - 1} (\pi \sqrt{D' D} / (c |L|)) & \text{if}\ \ D D' > 0\text{,} \\
  I_{2 s - 1} (\pi \sqrt{D' D} / (c |L|)) & \text{if}\ \ D D' < 0\text{.}
  \end{cases}
\end{align*}
The Kloosterman sum $K_{c,L}$ is defined in~(\ref{eq:klostermansum_definition}), and the $W$-Whittaker function $\mathcal{W}_{s, k - N/2}$ is given in~(\ref{eq:whittakerW_definition}).
\end{theorem}

Before proving this expansion, we use it to prove the Zagier-type duality announced in the introduction.
\begin{corollary}
\label{cor:zagiertypeduality}
The Fourier coefficients $c_{k,L,s}^{(n,r)} (n',r')$ of the Poincar\'e series $P_{k,L,s}^{(n,r)}$ of Proposition~\ref{prop:fourierexpansionofmaassjacobiforms} satisfy a Zagier-type duality with dual weights $k$ and $N + 2 - k$.  More precisely, suppose that $\Re(s) > 1 + N \slashdiv 2$, in which case the Poincar\'e series converges.  Then if $D, D' < 0$, there is a constant $h_{k,s}$ depending only on $(k, s)$ such that
\begin{align*}
  c_{k, L, s}^{(n,r)} (n',r')
& =
  h_{k,s} \, c_{N + 2 - k, L, s}^{(n', r')} (n, r)
\text{,}
\end{align*}
while if $D < 0$, $D' > 0$, there is a constant $h'_{k,s}$ depending only on $(k, s)$ such that
\begin{align*}
  c_{k, L, s}^{(n,r)} (n',r')
& =
  h'_{k,s} \, c_{N + 2 - k, L, s}^{(n', r')} (n, r)
\text{.}
\end{align*}
\end{corollary}

\begin{proof}
The corollary follows immediately from Theorem~\ref{th:poincarefourierexpansions} if we show that $K_{c,L}(n,r,n',r') = K_{c,L} (n', r', n, r)$.  This can be seen by changing $\lambda$ to $-d \lambda$ in~(\ref{eq:klostermansum_definition}).
\end{proof}

\begin{corollary}
If $s = k/2 - N/4$ (respectively, $1 + N/4 - k/2$), then the Poincar\'e series $P_{k,L,s}^{(n,r)}$ converges for $k > 2 + N$ (respectively, $k < 0$).  In both cases it is a Maa\ss-Jacobi form in $\mathbb{J}_{k,L}^{\mathrm{semi}}$ with Fourier expansion
\begin{align}
  P_{k, L, s}^{(n,r)} (\tau, z)
& =
  \wtM_{s, k - N/2}(\pi D y \slashdiv |L|) \; e(-i D y \slashdiv 4 |L|) \;
  \theta_{k, L}^{(r)} (\tau, t) \; q^n
\\[6pt]&
\nonumber
  \hspace{-1em}
  + \sum_{n' \in \ZZ, r' \in \ZZ^N, D' = 0} \hspace{-1em}
  c^{k} (n', r') \; q^{n'} \zeta^{r'}
\\[6pt]&
\nonumber
  \hspace{-1em}
  + \sum_{n' \in \ZZ, r' \in \ZZ^N, D' \ne 0} \hspace{-1em}
  c^{k} (n', r') \;
  e(-i D' y \slashdiv 4 |L|) \, \wtW_{s,k - N/2}(\pi D' y / |L|) \; q^{n'} \zeta^{r'}
\text{.}
\end{align}
\end{corollary}

\begin{proof}[of Theorem~\ref{th:poincarefourierexpansions}]
As usual we can choose as a system of representatives of $\JacFpara{N} \backslash \JacF{N}$ the elements
\begin{align*}
\label{eq:jacobigroupsystemofrepresentatives}
\Bigl(\left(\begin{matrix} a & b \\ c & d\end{matrix}\right), (a \lambda, b \lambda)\Bigr),
\end{align*}
where $(c,d)$ runs through coprime integers, $\lambda$ varies in $\ZZ^N$, and $a, b \in \ZZ$ are fixed for each pair $(c,d)$ such that $a d - b c = 1$.

The $\theta$-series in~(\ref{eq:poincarefourierexpansion}) arises naturally as the contribution of the representatives with $c = 0$.  Hence we must compute the contribution of representatives with $c \ne 0$.  Since the calculation is similar for both signs of $c$, we focus on $c > 0$.  We proceed by calculating the integral occurring in the Fourier transform.  First we separate as many terms as possible.  Toward this aim, note
\begin{gather*}
  \frac{a \tau + b}{c \tau + d}
 =
  \frac{a}{c} - \frac{1}{c (c \tau + d)}
\text{,} \qquad
  \frac{z}{c \tau + d} + \lambda \frac{a \tau + b}{c \tau + d}
=
  \frac{z - \lambda \slashdiv c}{c \tau + d} + \frac{a \lambda}{c}
\text{,}
\\[6pt]
  \frac{a \tau + b}{c \tau + d}\, \frac{L[\lambda]}{2}
  + \frac{\lambda L z + z L \lambda}{c \tau + d}
  - \frac{c L[z] \slashdiv 2}{c \tau + d}
=
  \frac{- c L[z - \lambda \slashdiv c] \slashdiv 2}{c \tau + d} + \frac{a L[\lambda]}{2 c}
\text{.}
\end{gather*}
Splitting fractions into their fractional and integer parts yields
\begin{align*}
  \sum_{\substack{c \in \ZZ_{>0} \\ d (c)^\times \\ \lambda \in \ZZ^N \mathop{\mathrm{mod}} c \ZZ^N \\
                   \alpha \in \ZZ,\, \beta \in \ZZ^N}}
&
  c^{-k} (\tau + c \slashdiv d + \alpha)^{-k}\;
  e \Bigl( \frac{-L[z - \lambda \slashdiv c + \beta]}{\tau + d \slashdiv c + \alpha} + a L[\lambda] \slashdiv c \Bigr)
\\&
  \cdot \phi_{k,L,s}^{(n,r)}
   \Bigl(a \slashdiv c - \frac{1}{c^2 (\tau + d \slashdiv c + \alpha)},\;
   \frac{z - \lambda \slashdiv c - \beta}{c (\tau + d \slashdiv c + \alpha)}
   + a \lambda \slashdiv c \Bigr) 
\text{.}
\end{align*}

Fixing the fractional parts and subsuming them under $\tau$, we need only compute the Fourier transform of
\begin{align*}
  \sum_{a \in \ZZ,\, \beta \in \ZZ^N} \!\!(\tau + \alpha)^{-k}\;
  e \Bigl(\frac{-L[z - \beta]}{\tau + \alpha}\Bigr)\;
  \phi_{k,L,s}^{(n,r)}
   \Bigl(a \slashdiv c - \frac{1}{c^2 (\tau + \alpha)},\;
   \frac{z - \beta}{c (\tau + \alpha)}
   + a \lambda \slashdiv c \Bigr) 
\text{.}
\end{align*}

By the Poisson summation formula for $L'$, the Fourier coefficients are
\begin{align*}
  a_y(n',r')
& =
  \int_{\RR \times \RR^N} t^{-k}\, e(L[w] \slashdiv t)\;\;
  \phi_{k,L,s}^{(n,r)} (a \slashdiv c - 1 \slashdiv c^2 t,\;
                     w \slashdiv c t + a \lambda \slashdiv c) \;
\\ &
  \hphantom{\int_{\RR \times \RR^N}} \qquad \cdot e(- n' x' + r' w) \; dx' du'
\text{,}
\end{align*}
where we use the variables $w = u' + i v' \in \CC^N$ and $t = x' + i y' \in \CC$.

Separating the real and imaginary part of
\begin{align*}
  a \slashdiv c - 1 \slashdiv c^2 t
& = 
  a \slashdiv c - x' \slashdiv c^2 |t|^2 + i(y \slashdiv c^2 |t|^2),
\end{align*}
we find that the integral equals
\begin{align*}
&
  e(n a \slashdiv c + a r \lambda \slashdiv c)
  \int_{\RR} t^{-k}\;
  \wtM_{s, k - N \slashdiv 2}(\pi D y \slashdiv |L| c^2 |t|^2)
\\&
  \qquad
  \cdot e(-n' x' + i L^{-1}[r] y\, c^2\, |t|^2
          - n x' \slashdiv c^2 |t|^2)
\\[4pt]&
  \qquad
  \cdot \int_{\RR^N} e(-r' w - L[w] \slashdiv t + r w \slashdiv c t ) \; du' dx'
\text{.}
\end{align*}
Evaluating the inner integral leaves us with
\begin{align*}
&
  \int_{\RR} t^{-(k - N/2)} \;
  \wtM_{s, k - N \slashdiv 2}(\pi D y \slashdiv |L|\, c^2\, |t|^2)
\\ &
  \hphantom{\int_{\RR}} \qquad
  \cdot e(D' x' \slashdiv 4 |L| - D x' \slashdiv 4 |L|\, c^2\, |t|^2) \; dx'
\text{,}
\end{align*}
which can be evaluated using \cite[page~176]{Fa77}.
\end{proof}
\vspace*{1.5em}

As a tool for the next proposition we will need skew-holomorphic Poincar\'e series.  For $k \ge 3$, set
\begin{align}
\label{eq:skewholomorphicpoincareseries_definition}
  P_{k,L} ^{(n,r),\mathrm{sk}}
:=
  \sum_{A \in \JacFpara{N} \backslash \JacF{N}}
  \! \left. e_{n,r,L}\, \right|_{1/2, k - 1/2, L} A,
\end{align}
where
\begin{align*}
  e_{n,r,L}(\tau, z)
:=
  e(n \tau + r z) \; e(-i D y \slashdiv 2 |L|)
\text{.}
\end{align*}

\begin{theorem}
The Poincar\'e series $P_{k,L} ^{(n,r),\mathrm{sk}}$ defined in~(\ref{eq:skewholomorphicpoincareseries_definition}) is an element of\/ $\mathbb{J}_{k,L}^{\mathrm{sk},\mathrm{cusp}}$ and has Fourier expansion
\begin{align*}
  P_{k,L} ^{(n,r),\mathrm{sk}}(\tau, z)
& =
  e(-i D y \slashdiv 2 |L|) \; \theta_{k-1, L}^{(r)} (\tau, z) \; q^n
\\ &
  + \sum_{\substack{n' \in \ZZ, r' \in \ZZ^N \\ D > 0}}
    \!\! c(n',r') \; e(-i D' y \slashdiv 2 |L|) \; q^{n'} \zeta^{r'}
\text{.}
\end{align*}
Here the $\theta$-series is defined in~(\ref{eq:thetaseries_definition}), and the coefficients $c(n',r')$ are
\begin{align*}
  c(n', r')
& = 
  b(n', r') + (-1)^k b(n', -r'),
\end{align*}
where $b$ depends on $D$ and $D'$.  We have
\begin{align*}
  b(n', r')
& =
  2^{1 - N / 2} \pi i^{-k + 1}\, |L|^{-1 \slashdiv 2} \;
  (D' \slashdiv D)^{k / 2 - (N + 2) / 4}
\\ &
  \sum_{c \in \ZZ_{>0}} \! c^{-(N + 2) / 2} \; K_{c,L}(n, r, n', -r') \;
  J_{k - (N + 2) / 2} (\pi \sqrt{D D'} \slashdiv |L| c)
\text{,}
\end{align*}
the Kloosterman sum $K_{c,L}$ being defined in~(\ref{eq:klostermansum_definition}).
\end{theorem}

\begin{proof}
The proof is analogous to that of Theorem~\ref{th:poincarefourierexpansions}.
\end{proof}

\begin{remark}
The Poincar\'e series $P_{k,L}^{(n,r),\mathrm{sk}}$ span\/ $\mathbb{J}_{k,L}^{\mathrm{sk},\mathrm{cusp}}$.  This can be seen as in \cite{Sk90}: evaluation of the Petersson scalar product of an arbitrary form $f$ with $P_{k,L}^{(n,r),\mathrm{sk}}$ yields the $(n,r)^\thup$ Fourier coefficient of $f$.  Hence any cusp form orthogonal with respect to the Petersson scalar product to all Poincar\'e series vanishes.
\end{remark}

\begin{proposition}
\label{prop:poincaremapstopoincare}
The Maa\ss-Jacobi Poincar\'e series $P_{k,L,k \slashdiv 2 - N \slashdiv 4} ^{(n,r)}$ with $k > 2 + N $ is meromorphic.

For $k < 0$, $\xi_L$ maps the Maa\ss-Jacobi Poincar\'e series $P_{k,L,1 + N \slashdiv 4 - k \slashdiv 2}^{(n,r)}$ to the skew-holomorphic Poincar\'e series $P_{2 + N - k,L}^{(n,r),\mathrm{sk}}$ up to a constant factor.
\end{proposition}

\begin{proof}
This follows immediately if we recall that
\begin{align*}
  \wtW_{1 + N \slashdiv 4 - k \slashdiv 2, k - N \slashdiv 2} (y)
& =
  \begin{cases}
    e^{-y \slashdiv 2} &\text{ if }y > 0, \\
    e^{-y \slashdiv 2} \Gamma((N + 2) \slashdiv 2 - k, -y) &\text{ otherwise.}
  \end{cases}
\end{align*}
\end{proof}

We obtain the following simple but important corollary.
\begin{corollary}
\label{cor:xisurjectivityandimages}
For $k < 0$ the restriction of the $\xi$-operator to $\mathbb{J}_{2 + N - k,L}^{\mathrm{semi},\mathrm{cusp}}$ is surjective onto $\mathbb{J}_{k,L}^{\mathrm{sk},\mathrm{cusp}}$.
\end{corollary}
\begin{proof}
Applying the $\xi$-operator and exploiting the fact that skew-holomorphic Poincar\'e series span $\mathbb{J}_{k,L}^{\mathrm{sk},\mathrm{cusp}}$, the result follows from Proposition~\ref{prop:poincaremapstopoincare}.
\end{proof}


%% file: casimiroperators.tex
\def\cS{{\cal S}}
\def\cdo{covariant \dog}

\def\ad{\mathop{\rm ad}\nolimits}
\def\CDO{\mathop{\rm CDO}\nolimits}
\def\End{\mathop{\rm End}\nolimits}
\def\Frac{\mathop{\rm Frac}\nolimits}
\def\IDO{\mathop{\rm IDO}\nolimits}
\def\ie{{\em i.e.,\/}}
\def\Sym{\mathop{\rm Sym}\nolimits}

\section{The Casimir operator} \label{sec:casimiroperators}

The purpose of this section is to prove Theorem~\ref{thm:actionofcenter} and Propositions~\ref{prop:RLops}, \ref{prop:allCDOs}, and~\ref{prop:allIDOs}.  Recall the real Jacobi group $\t G^J_N$ defined in Section~\ref{sec:maassjacobiforms}, and let $\t\dg^J_N$ be its complexified Lie algebra.  Making identifications as in Definition~\ref{def:realjacobicentralextension}, we find
\begin{gather*}
   \t\dg^J_N\, :=\, \bigl\{ (M, X, \kappa):\,
   M \in \dsl_2(\bC),\, X \in \Mat{N,2}{\bC},\,
   \kappa \in \MatT{N}{\bC} \bigr\}
\text{,}
\\[6pt]
   \bigl[ (M, X, \kappa),\, (\check M, \check X, \check \kappa) \bigr]\, =\, 
   \bigl( [M, \check M],\, X \check M - \check X M,\,
   \check X J_2 X^\T - X J_2 \check X^\T \bigr)
\text{.}
\end{gather*}

The exponential map $\exp: \t\dg^J_N \to \t G^J_N$ is
\begin{equation} \label{eq:exp}
   \exp(M, X, \kappa)\, =\, (e^M,\, X g(M),\, \kappa - X h(M) J_2 X^\T),
\end{equation}
where $g(z) := (e^z - 1)/z$ and $h(z) := (e^z - z - 1)/z^2$.  This formula can be proven efficiently by checking that $\t\dg^J_N$ embeds in $\dgl_{2N+2}(\bR)$ as follows.  For reference, we also give the embedding of $\t G^J_N$ into $\GL{2N+2}(\bR)$:
\begin{equation*}
   (M, X, \kappa)_{\t\dg^J_N} \mapsto
   \left( \begin{matrix} 0 & X & \kappa \\
                           & M & -J_2 X^\T \\ & & 0 \end{matrix} \right),
   \quad
   (M, X, \kappa)_{\t G^J_N} \mapsto
   \left( \begin{matrix} I_N & X & \kappa \\
                           & M & -M J_2 X^\T \\ & & I_N \end{matrix} \right).
\end{equation*}
Observe that conjugation by the matrix of the cyclic permutation $(1, 2, \ldots, N+1)$ carries $\t G^J_N$ into a minimal parabolic subgroup of the standard realization of $\Sp_{2N+2}(\bR)$, the image being everything but the center.

We note that the formula for the exponential function given in Section~5.2 of \cite{BCR07} is incorrect.  It should be identical to~(\ref{eq:exp}), but the $h(M)$ term was missed.  (It is given correctly in Section~5.1 of \cite{BCR10}.)  The embedding $\pi_4$ on p.~148 of \cite{BCR07} is also wrong: the $(1,4)$ and $(3,4)$ entries should make up $-M J_2 X^\T$ rather than $-J_2 X^\T$.

In order to fix a basis for $\t\dg^J_N$, recall that we write $\epsilon_{ij}$ for the elementary matrix with $(i, j)^\thup$ entry~$1$ and other entries~$0$, the size of the matrix being determined by the context.  The standard basis of $\dsl_2$ is of course $E := \epsilon_{12}$, $F := \epsilon_{21}$, and $H := \epsilon_{11} - \epsilon_{22}$.  For a basis of $\rmM_{N,2}$ we take $e_i := \epsilon_{i2}$ and $f_i := \epsilon_{i1}$, and for $\rmM^\T_N$ we take $Z_{ij} := \oh (\epsilon_{ij} + \epsilon_{ji})$.  Then $\t\dg^J_N$ has basis
$$ \bigl\{ E, F, H;\ e_i, f_i,\, 1 \le i \le N;\
           Z_{ij},\, 1 \le i \le j \le N \bigr\}. $$
The brackets of this basis are as follows.  Those on $\dsl_2$ are standard, and $\rmM^\T_N$ is the center $\dz(\t\dg^J_N)$.  Under $\ad(H)$, the $e_i$ are of weight~$1$ and the $f_i$ are of weight~$-1$.  The $\ad(E)$ and $\ad(F)$ actions are given by
$$ \ad(E):\ e_i \mapsto 0,\, f_i \mapsto -e_i; \quad
   \ad(F):\ e_i \mapsto -f_i,\, f_i \mapsto 0. $$
Finally, $[e_i, f_j] = -2 Z_{ij}$.

The first step in proving Theorem~\ref{thm:actionofcenter} is to compute the center $\dZ(\t\dg^J_N)$ of the \uea\ $\dU(\t\dg^J_N)$.  Towards this end, let us write $e$ and $f$ for the column vectors with entries $e_i$ and $f_i$, respectively, and $Z$ for the symmetric matrix with entries $Z_{ij}$.  This permits us to write conveniently such elements of $\dU(\t\dg^J_N)$ as $e^\T Z f$ and $\det(Z)$, the determinant of $Z$.

Observe that $\det(Z) e^\T Z^{-1} f$ is a well-defined element of $\dU(\t\dg^J_N)$.  We will need in addition the following more subtle fact:
$$ P\, :=\, \det(Z) \bigl( e^\T (e^\T Z^{-1} f) Z^{-1} f
   - (e^\T Z^{-1} e) (f^\T Z^{-1} f) \bigr) $$
is an element of $\dU(\t\dg^J_N)$ of degree $N + 2$.  To prove this, note that $\det(Z) P$ is clearly in $\dU(\t\dg^J_N)$.  Check that if we specialize $Z$ to a diagonal matrix of scalars, then $\det(Z) = 0$ implies $\det(Z) P = 0$.  Since the symmetric determinant is an irreducible polynomial, the result follows.

We now define the {\em Casimir element\/} of $\dU(\t\dg^J_N)$, which has degree $N + 2$:
\begin{equation*} \begin{array}{rl}
   \Omega_N\ :=\ \det(Z) \hskip-10pt
   & \Bigl( H^2 - (N + 2)H + 4EF \\[6pt]
   &\ -\, \bigl( H - \oh(N + 3) \bigr) e^\T Z^{-1} f
   + E f^\T Z^{-1} f - e^\T Z^{-1} e F \\[6pt]
   &\ +\, {\ts\frac{1}{4}} e^\T (e^\T Z^{-1} f) Z^{-1} f
   - {\ts\frac{1}{4}} (e^\T Z^{-1} e) (f^\T Z^{-1} f) \Bigr).
\end{array} \end{equation*}

\begin{theorem} \label{thm:center}
The center $\dZ(\t\dg^J_N)$ is polynomial on $1 + {N + 1 \choose 2}$ generators:
$$ \dZ(\t\dg^J_N)\, =\, \bC[\Omega_N, Z_{ij}:\, 1 \le i \le j \le N]. $$
\end{theorem}

In order to prove this theorem we will describe {\em virtual copies of semisimple subalgebras,\/} which provide a simple but very clever way to compute the centers of the \uea s of a certain class of semidirect sum Lie algebras.  This idea was discovered by Borho \cite{Bo76} for $\t\dg^J_1$, and independently by Quesne \cite{Qu88} for various Lie algebras of 1-dimensional center.  It was generalized by Campoamor-Stursburg and Low \cite{CL09} and applied to further examples.

Before we give the proof of Theorem~\ref{thm:center} arising from the virtual copy of $\dsl_2$, let us outline a straightforward but considerably less efficient proof.  For any Lie algebra $\dg$, the {\em symmetrizer map\/} $\Sym$ from the symmetric algebra $\cS(\dg)$ to $\dU(\dg)$ restricts to a vector space isomorphism from the invariants $\cS(\dg)^\dg$ to $\dZ(\dg)$.  The strategy is to compute $\cS(\t\dg^J_N)^{\dsl_2}$ and then locate $\cS(\t\dg^J_N)^{\t\dg^J_N}$ inside it.  

Define the following elements of $\cS(\t\dg^J_N)$: $Q_0 := H^2 + 4FE$, and
\begin{equation*}
   Q_{ij}\, :=\, \oh (f_i e_j - f_j e_i), \quad
   C_{ij}\, :=\, E f_i f_j - \oh H (f_i e_j + f_j e_i) - F e_i e_j
\end{equation*}
for $1 \le i, j \le N$.  It can be shown that together with the $Z_{ij}$, these elements generate $\cS(\t\dg^J_N)^{\dsl_2}$.  For reference, we note that they satisfy the relations
\begin{equation*}
   Q_{ij} Q_{kl} + Q_{il} Q_{jk} + Q_{ik} Q_{lj}\, =\, 0, \quad
   C_{ij} C_{kl} - C_{ik} C_{jl} + Q_0 Q_{il} Q_{jk}\, =\, 0
\end{equation*}
for $1 \le i, j, k, l \le N$.  Thus $\cS(\t\dg^J_N)^{\t\dg^J_N}$ is the kernel of the actions of the $\ad(e_i)$ and $\ad(f_i)$ on the algebra generated by $\bigl\{ Q_0, Q_{ij}, C_{ij}, Z_{ij} \bigr\}_{ij}$. 

Suppose that $P$ is an element of this kernel.  Regard it as a polynomial in the $Z_{ij}$ with coefficients generated by $\bigl\{ Q_0, Q_{ij}, C_{ij} \bigr\}_{ij}$.  One can use the $\ad(e_i)$- and $\ad(f_i)$-invariance of $P$ to prove that the coefficients of the monomials of top $Z$-degree lie in $\bC[Q_0]$.

Now let $Q$ and $C$ be the matrices with entries $Q_{ij}$ and $C_{ij}$, and define
\begin{equation*}
   P_N\, :=\, \det(Z) \bigl( Q_0 + \tr (Z^{-1} C) + \oh \tr (Z^{-1} Q)^2 \bigr).
\end{equation*}
This is polynomial in $Z$ because $Q$ is skew-symmetric, and direct verification shows that it is $\t\dg^J_N$-invariant.

Collecting these results, one finds easily that $P_N$ generates $\cS(\t\dg^J_N)^{\t\dg^J_N}$ over the extension of $\bC[Z_{ij}]$ by $\det(Z)^{-1}$.  With a little more work one can drop the necessity to adjoin $\det(Z)^{-1}$.

The last step is to symmetrize $P_N$.  A long computation leads to $\Sym(P_N) = \Omega_N + \frac{1}{4} N(N+3) \det(Z)$, proving Theorem~\ref{thm:center}.

\subsection{Virtual semisimple subalgebras} \label{ssec:VSSs}

This section is an exposition of some of the results of \cite{Bo76, Qu88, CL09}.  Let $\dg$ be any complex finite dimensional Lie algebra.  Write $\ds \oplus_s \dr$ for its semidirect sum Levi decomposition, in which its semisimple part $\ds$ acts on its solvable part $\dr$.  Let $\dz$ be the center $\dz(\dg)$, which is of course contained in $\dr$.  For any subalgebra $\dh$ of $\dg$ containing $\dz$, define
$$ \dU_\dz(\dh)\, :=\, \dU(\dh) \otimes_{\dU(\dz)} \Frac\bigl(\dU(\dz)\bigr), \quad
   \dZ_\dz(\dh)\, :=\, \dZ(\dh) \otimes_{\dU(\dz)} \Frac\bigl(\dU(\dz)\bigr). $$

Throughout this section we make the following assumption:
\begin{equation} \label{eq:ansatz}
   \mbox{\rm There is a Lie algebra homomorphism $\eta: \dg \to \dU_\dz(\dr)$
             with $\eta|_\dr = 1$.}
\end{equation}
We will see that under this assumption, there is a {\em virtual copy\/} $\ds_\nu$ of $\ds$ in $\dU_\dz(\dg)$ such that $\dZ(\ds_\nu) \otimes \dZ_\dz(\dr)$ and $\dZ_\dz(\dg)$ are equal {\em as algebras.\/}  This greatly simplifies the deduction of $\dZ(\dg)$.

To our knowledge, general conditions on $\dg$ under which $\eta$ exists are not known.  Such conditions would be interesting.  The case that $\dg$ is perfect and maximally centrally extended, \ie\ $[\dg, \dg] = \dg$ and $H^2(\dg, \bC) = 0$, may be significant.  However, although $\t\dg^J_N$ has these properties, $\eta$ does not exist for all such Lie algebras.  The semidirect product of $\dsl_2$ with its adjoint \rep\ is a counterexample.

\begin{lemma} \label{lemma:virtual s}
Assuming~(\ref{eq:ansatz}), $\eta$ commutes with the $\ad$-action of $\dg$.
\end{lemma}

\begin{proof}
First prove that $\eta$ is an $\dr$-map, and note that $\eta|_\dr = 1$ is an $\ds$-map.  Then prove that $[S, \Theta] = [\eta(S), \Theta]$ for all $S \in \ds$ and $\Theta \in \dU_\dz(\dr)$.  Finally, prove that $\eta|_\ds$ is an $\ds$-map.
\end{proof}

The proofs of the next two results are left to the reader.

\begin{corollary}
Assuming~(\ref{eq:ansatz}), the map $\nu := 1 - \eta: \dg \to \dU_\dz(\dg)$ is a $\dg$-map and a Lie algebra homomorphism.  It annihilates $\dr$ and is injective on $\ds$.  Its image $\ds_\nu := \nu(\ds)$, the\/ {\em virtual copy\/} of $\ds$, commutes with $\dr$.  The natural multiplication map defines algebra isomorphisms
$$ \dU(\ds_\nu) \otimes \dU_\dz(\dr)\, \cong\, \dU_\dz(\dg), \quad
   \dZ(\ds_\nu) \otimes \dZ_\dz(\dr)\, \cong\, \dZ_\dz(\dg). $$
\end{corollary}

\begin{lemma}
Assumption~(\ref{eq:ansatz}) holds for $\t\dg^J_N$: define $\eta$ by
\begin{eqnarray*}
   \eta(E) &:=& {\ts\frac{1}{4}} e^\T Z^{-1} e, \qquad
   \eta(F)\ :=\ -{\ts\frac{1}{4}} f^\T Z^{-1} f, \\[6pt]
   \eta(H) &:=& {\ts\frac{1}{4}} (e^\T Z^{-1} f + f^\T Z^{-1} e)
           \ =\ {\ts\frac{1}{2}} (N + e^\T Z^{-1} f).
\end{eqnarray*}
\end{lemma}

\medbreak \noindent {\bf Proof of Theorem~\ref{thm:center}.}
Take $\dg = \t\dg^J_N$, where $\ds = \dsl_2$.  As is well known, the Casimir operator $\Omega_{\dsl_2} := H^2 - 2H + 4EF$ generates $\dZ(\ds)$, and so $\nu(\Omega_{\dsl_2}) = \nu(H)^2 - 2\nu(H) + 4\nu(E) \nu(F)$ generates $\dZ(\ds_\nu)$.  A short calculation gives $\nu(\Omega_{\dsl_2}) = \det(Z)^{-1} \Omega_N + \frac{1}{4} N (N + 4)$.  It is not hard to verify that $\dZ(\dr)$ is $\dU(\dz)$.  Clearing denominators and noting that we have polynomial independence, the result follows.  \hfill $\Box$

\subsection{The slash actions of $\t\dg^J_N$} \label{ssec:central action}

Recall from Section~\ref{sec:maassjacobiforms} that associated to each scalar cocycle $\gamma: \tJac{N} \to C^\infty(\HS_{1,N})$ there is a right action $|_\gamma$ of $\tJac{N}$ on $C^\infty(\HS_{1,N})$, defined for $g \in \tJac{N}$ and $x \in \HS_{1,N}$ by $\phi|_\gamma [g](x) := \gamma(g, x) \phi(g x)$.

We will denote the differential of $\gamma$ at the identity, a linear map from $\t\dg^J_N$ to $C^\infty(\HS_{1,N})$, by the same symbol $\gamma$.  Thus $\gamma(Y, x) := \partial_t|_{t=0} \gamma(e^{tY}, x)$ for $Y \in \t\dg^J_N$.  The differential right action $|_\gamma$ of $\t\dg^J_N$ on $C^\infty(\HS_{1,N})$ is
$$ \phi|_\gamma [Y](x)\, :=\, 
   \partial_t|_{t=0}\, \bigl( \gamma(e^{tY}, x)\, \phi(e^{tY} x) \bigr)\, =\,
   \gamma(Y, x)\, \phi(x) + \partial_t|_{t=0}\, \phi(e^{tY} x). $$
This action extends to $\dU(\t\dg^J_N)$ as usual, and elements of $\dU(\t\dg^J_N)$ of order~$r$ act by \dog s of order~$\le r$.

Following Definition~\ref{defn:CDOs}, let $\bD_\gamma$ denote the algebra of \dog s on $\HS_{1,N}$ invariant with respect to the action $|_\gamma$ of $\tJac{N}$.  Since $\tJac{N}$ is connected, $\bD_\gamma$ is precisely those operators commuting with the $|_\gamma$-action of $\t\dg^J_N$.  In particular, $|_\gamma$ maps the center $\dZ(\t\dg^J_N)$ of $\dU(\t\dg^J_N)$ into the center $Z(\bD_\gamma)$ of $\bD_\gamma$.  (In \cite{BCR10} it is proven that $\dZ(\t\dg^J_1)$ covers $Z(\bD_\gamma)$ for all scalar cocycles of $\tJac{1}$.  It would be interesting to decide this question for $N>1$.)

\medbreak \noindent {\bf Proof of Theorem~\ref{thm:actionofcenter}.}
Recall the cocycles $\beta^k \alpha_L$ from Lemma~\ref{lemma:cocycles} defining the slash actions $|_{k,L}$.  Easy calculations using~(\ref{eq:exp}) give the actions of our basis of $\t\dg^J_N$: writing $\sL$ for $2 \pi i L$ and $A[w]$ for $w^\T A w$ as earlier,
\begin{equation} \label{eq:Lie slash action}
\begin{array}{rcl}
   |_{k,L} [E] &=& 2 \Re(\partial_\tau), \\[6pt]
   |_{k,L} [F] &=& -2 \Re \tau (\tau \partial_\tau + z^\T \partial_z)
                   - k \tau - \sL[z], \\[6pt]
   |_{k,L} [H] &=& 2 \Re(2\tau \partial_\tau + z^\T \partial_z) + k, \\[6pt]
   |_{k,L} [0, X, \kappa] &=& 2\Re (X_1 \tau + X_2)^\T \partial_z
                              + 2 X_1^\T \sL z + \tr(\kappa \sL).
\end{array}
\end{equation}
In using these formulas to compute the $|_{k,L}$-actions of elements of $\dU(\t\dg^J_N)$, care must be taken to reverse the order of multiplication, because $|_{k,L}$ is a right action but the formulas~(\ref{eq:Lie slash action}) are given in left-acting notation.  Thus for example $\phi|[Y_1 Y_2]$ is $\bigl(\phi|[Y_1]\bigr)|[Y_2]$.

Since $|_{k,L}[Z_{ij}] = \sL_{ij}$, Theorem~\ref{thm:center} shows that $|_{k,L}[\Omega_N]$ generates the $|_{k,L}$-action of $\dZ(\t\dg^J_N)$.  Using~(\ref{eq:Lie slash action}), a straightforward but long computation gives
\begin{equation*}
   |_{k,L}[\Omega_N]\, =\, \det(\sL) \bigl( k (k - N - 2) - 2\cC^{k,L} \bigr). \qquad \Box
\end{equation*}

\subsection{Covariant \dog s} \label{ssec:CDOs}
In order to prove Propositions~\ref{prop:RLops}, \ref{prop:allCDOs}, and~\ref{prop:allIDOs}, we must recall the algebraic side of the general theory of \ido s (IDOs) developed by Helgason in the 1950's (see, {\em e.g.,\/} Section~2 of \cite{Helg77}).  Here we will adapt the framework of Section~4 of \cite{BCR10} from IDOs to {\em \cdo s\/} (CDOs), by regarding them as nilpotent IDOs on the direct sum of the range and domain spaces.  Thus to treat scalar-valued CDOs we must consider vector-valued IDOs.

Let $G$ be a real Lie group, $K$ a closed subgroup, and $V$ a complex vector space.  Given $x \in G$, denote the coset $xK$ by $\ov x$.  A {\em $V$-valued 1-cocycle of\/ $G$ on\/ $G/K$\/} is a smooth function
$$ \gamma: G \times G/K \to \GL{}(V)
   \mbox{\rm\ satisfying\ }
   \gamma(gh, \ov x) = \gamma(h, \ov x) \gamma(g, h \ov x) $$
for all $g,\, h,\, x\, \in G$.  The associated right action of $G$ on $C^\infty(G/K) \otimes V$ is
$$ f|_\gamma[g](\ov x)\, :=\, \gamma(g, \ov x) f(g \ov x), $$
and the associated \rep\ of $K$ on $V$ is $\pi_\gamma(k) := \gamma(k, \ov e)^{-1}$.

Suppose that $V'$ is a vector space of the same dimension as $V$, and $\gamma'$ is a $V'$-valued 1-cocycle of $G$ on $G/K$.  Then $\gamma'$ is said to be {\em cohomologous\/} to $\gamma$ if there is a smooth map $b$ from $G/K$ to the set of invertible linear maps from $V$ to $V'$ such that $\gamma'(g, \ov x) b(g \ov x) = b(\ov x) \gamma(g, \ov x)$.  In this case $f \mapsto b f$ is an equivalence from $|_\gamma$ to $|_{\gamma'}$, and $b(\ov e)$ is an equivalence from $\pi_\gamma$ to $\pi_{\gamma'}$.  Conversely, if $\pi_\gamma$ and $\pi_{\gamma'}$ are equivalent, then $\gamma$ and $\gamma'$ are cohomologous.  If $G/K$ is simply connected, then given any complex finite dimensional \rep\ $\pi$ of $K$ there exists a cocycle $\gamma$ such that $\pi_\gamma = \pi$, and so one has a natural bijection between slash actions of $G$ on $G/K$ and \rep s of $K$.

We now define CDOs in the general setting.  Let $V$ and $V'$ be any two vector spaces, not necessarily related.  Fix 1-cocycles $\gamma$ and $\gamma'$ of $G$ on $G/K$ taking values in $V$ and $V'$, respectively.

\begin{definition} \label{defn:generalCDOs}
A \dog\/ $T:C^\infty(G/K) \otimes V \to C^\infty(G/K) \otimes V'$ is\/ {\em covariant from\/ $|_\gamma$ to\/ $|_{\gamma'}$\/} if for all\/ $g \in G$ and\/ $f \in C^\infty(G/K) \otimes V$, we have
$$ T \bigl( f|_\gamma [g] \bigr)\ =\ (T f) \big|_{\gamma'} [g]. $$

Let\/ $\bD_{\gamma, \gamma'}(G/K)$ be the space of CDOs from\/ $|_\gamma$ to\/ $|_{\gamma'}$, and let\/ $\bD^r_{\gamma, \gamma'}(G/K)$ be the space of those of order~$\le r$.  When $\gamma = \gamma'$, we refer to such operators as\/ {\em $|_\gamma$-invariant\/} and write simply\/ $\bD_\gamma(G/K)$ and\/ $\bD^r_\gamma(G/K)$.
\end{definition}

The following proposition adapts Section~4.3 of \cite{BCR10} to CDOs.  Let $\dg$ and $\dk$ be the complexified Lie algebras of $G$ and $K$, and assume that the pair $K \subseteq G$ is {\em reductive,\/} \ie\ there exists a $K$-splitting $\dk \oplus \dm$ of $\dg$.  Recall that $\cS$ denotes the symmetric algebra and superscripts indicate invariants.

\begin{proposition} \label{prop:generalCDOs}
There exists a filtration-preserving linear bijection
$$ \CDO_{\gamma, \gamma'}: \bigl(\cS(\dm) \otimes \Hom(V, V') \bigr)^K
   \to \bD_{\gamma, \gamma'}(G/K). $$
It is compatible with multiplication at the symbol level: if $\gamma''$ is a third cocycle taking values in a space $V''$ and $\Theta$ and $\Theta'$ are $K$-invariant elements of $\cS(\dm) \otimes \Hom(V, V')$ and $\cS(\dm) \otimes \Hom(V', V'')$, respectively, then $\CDO_{\gamma, \gamma''}(\Theta' \Theta)$ and $\CDO_{\gamma', \gamma''}(\Theta') \circ \CDO_{\gamma, \gamma'}(\Theta)$ have the same symbol.
\end{proposition}

\begin{proof}
Construct a $V \oplus V'$-valued cocycle $\gamma \oplus \gamma'$ in the obvious way.  Equation~(38) of \cite{BCR10} defines a filtration-preserving linear bijection
$$ \IDO_{\gamma \oplus \gamma'}: \bigl(\cS(\dm) \otimes \End(V \oplus V') \bigr)^K
   \to \bD_{\gamma \oplus \gamma'}(G/K) $$
which is an algebra isomorphism at the symbol level.  Tracing the definitions leading to~(38) shows that $\IDO_{\gamma \oplus \gamma'}$ restricts to the desired map $\CDO_{\gamma, \gamma'}$.
\end{proof}

We now give a general result (probably already known) of independent interest: in the reductive case, the CDOs of order~1 generate all CDOs.  We do not know if it holds in the absence of reductivity.

\begin{proposition} \label{prop:ord1CDOs}
Let $K \subseteq G$ be reductive.  Then all CDOs of order~$r$ are linear combinations of compositions of up to~$r$ CDOs of order~1.
\end{proposition}

\begin{proof}
By Proposition~\ref{prop:generalCDOs}, it suffices to show that $\bigl(\cS^r(\dm) \otimes \Hom(V, V') \bigr)^K$ is contained in the product
$$ \Hom(V_r, V')^K
   \bigl(\dm \otimes \Hom(V_{r-1}, V_r) \bigr)^K \cdots
   \bigl(\dm \otimes \Hom(V_1, V_2) \bigr)^K
   \bigl(\dm \otimes \Hom(V, V_1) \bigr)^K $$
for some \rep s $V_1, \ldots, V_r$ of $K$ (the first factor contains order~0 operators and can be merged with the second factor).

Set $V_s := \cS^s(\dm^*) \otimes V$.  Fix any basis $\{X_j\}$ of $\dm$, and let $\{X^*_j\}$ be the dual basis of $\dm^*$.  Let $I_s$ be $\sum_j X_j X^*_j$, regarded as an element of $\dm \otimes \Hom(V_{s-1}, V_s)$ in the natural way.  Verify that as such, it is $K$-invariant.  We will in fact prove that $\bigl(\cS^r(\dm) \otimes \Hom(V, V') \bigr)^K$ is equal to $\Hom(V_r, V')^K I_r \cdots I_1$.

For this, it is enough to show that right composition with $I_r \cdots I_1$ is an injection from $\Hom(V_r, V')$ to $\cS^r(\dm) \otimes \Hom(V, V')$, as these two \rep s of $K$ are equivalent.  Using (subscript) monomial notation $X_J$ and $X^*_J$, observe that $I_r \cdots I_1$ may be written as $\sum_J {r \choose J} X_J X^*_J$.  Given $H$ in $\Hom(V_r, V')$, regard it as a map from $\cS^r(\dm^*)$ to $\Hom(V, V')$.  Then $H I_r \cdots I_1$ is $\sum_J {r \choose J} X_J \otimes H(X^*_J)$, proving the injectivity.
\end{proof}

The construction in the preceding proof is inefficient in practice.  Let us describe a more useful approach in the case that $K$ is abelian, $G/K$ is simply connected, and we restrict to cocycles $\gamma$ such that $\pi_\gamma$ is a completely reducible \rep\ of $K$.  Here the irreducible \rep s of $K$ are 1-dimensional, so it suffices to prove the result for CDOs between scalar slash actions.  Given a scalar 1-cocycle $\gamma$, write $\bC_\gamma$ for $\bC$ endowed with the $K$-action $\pi_\gamma$, and $1_\gamma$ for $1 \in \bC_\gamma$.  If $\gamma'$ is a second cocycle, then $\gamma'/\gamma$ is again a cocycle and $\Hom(\bC_\gamma, \bC_{\gamma'})$ is $K$-equivalent to $\bC_{\gamma'/\gamma}$.  Therefore by Proposition~\ref{prop:generalCDOs}, there is an order-preserving bijection
$$ \bD_{\gamma, \gamma'}(G/K) \cong
   \bigl( \cS(\dm) \otimes \bC_{\gamma'/\gamma} \bigr)^K. $$

The space on the right is essentially the $\pi_{\gamma/\gamma'}$-isotype of $\cS(\dm)$.  Let $\{X_j\}$ be a $K$-eigenbasis of $\dm$, and for each $j$ let $\chi_j$ be a scalar cocycle such that $\bC X_j$ is a copy of $\pi_{1/\chi_j}$ under $K$ (such $\chi_j$ exist because $G/K$ is simply connected).  Then for each scalar cocycle $\gamma$ we have the order~1 CDO
\begin{equation} \label{eq:CDObasis}
   X^\gamma_j\, :=\, \CDO_{\gamma, \gamma \chi_j} (X_j \otimes 1_{\chi_j})
   \,\in\, \bD^1_{\gamma, \gamma \chi_j}(G/K).
\end{equation}

Note that by Lemma~4.6 of \cite{BCR10}, the symbol of $X^\gamma_j$ is independent of $\gamma$.  For clarity we will sometimes use the notation $X^{\CDO}_j$ for the diagonal action of $\oplus_\gamma X^\gamma_j$ on the algebraic direct sum of all the scalar slash actions.  

Combining the preceding discussion with Propositions~\ref{prop:generalCDOs} and~\ref{prop:ord1CDOs} gives the following corollary (the last statement of which follows from an examination of the definition of $\CDO_{\gamma, \gamma'}$; see \cite{BCR10}).

\begin{corollary} \label{cor:CDObasis}
Assume that $K \subseteq G$ is reductive, $K$ is abelian, and $G/K$ is simply connected.  Let $\gamma$ and $\gamma'$ be scalar cocycles.  Then\/ 
\begin{equation*}
   \bD^r_{\gamma, \gamma'}(G/K) \mbox{\rm\ \ has basis\ \ }
   \ts \bigl\{ \prod_j (X^{\CDO}_j)^{J_j}:\ \prod_j \chi_j^{J_j} 
   = \gamma'/\gamma,\ \sum_j J_j \le r \bigr\}.
\end{equation*}

The symbol of $X^{\CDO}_j$ at the base point $\ov e$ coincides with that of the left action of $X_j$, which in turn coincides with that of\/ $-|_\gamma [X_j]$ for all $\gamma$.
\end{corollary}

We remark that the various $X^{\CDO}_j$ do not necessarily commute, but since we are in the scalar setting their commutators are of order~$\le 1$.  More precisely,
\begin{equation} \label{eq:commutators}
   X_{j'}^{\gamma \chi_j} \circ X_j^\gamma -
   X_j^{\gamma \chi_{j'}} \circ X_{j'}^\gamma\ \in\
   \Span \bigl\{ X^\gamma_{j''}:\, \chi_{j''} = \chi_{j'} \chi_j \bigr\}
   \, \oplus \, \delta_{1, \chi_{j'} \chi_j} \bC 1.
\end{equation}
(The Kronecker delta coefficient of the final summand $\bC 1$ indicates that it should be omitted unless $\chi_{j'} \chi_j =1$.)  In light of the basis given in Corollary~\ref{cor:CDObasis}, this leads to the following proposition.

\begin{proposition} \label{prop:CDOrelations}
For $K \subseteq G$ reductive, $K$ abelian, and $G/K$ simply connected, the ``reordering relations''~(\ref{eq:commutators}) generate all relations between CDOs.
\end{proposition}

We would be interested to know if the CDO relations of order~2 generate all CDO relations in more general settings, for example for $K$ non-abelian.

\subsection{Proofs of CDO propositions from Section~\ref{sec:maassjacobiforms}}
\label{ssec:CDOproofs}
We now specialize to the setting of Section~\ref{sec:maassjacobiforms}.  Recall that the subgroup of $\tJac{N}$ stabilizing the point $(i, 0)$ in $\HS_{1, N}$ is $\t K^J_N = \SO_2 \times \{0\} \times \MatT{N}{\bR}$.  This group has complex Lie algebra $\t\dk^J_N := \dgo_2 \times \{0\} \times \MatT{N}{\bC}$.

Define a linear map $\tau: \t\dg^J_N \to \t\dg^J_N$ by $\tau(X) = \t X$, where
\begin{equation*}  \begin{array}{lll}
   \t H := i(F-E), \quad
   &\t E := {\ts\frac{1}{2}} \bigl(H + i(F + E)\bigr), \quad
   &\t F := {\ts\frac{1}{2}} \bigl(H - i(F + E)\bigr), \\[6pt]
   \t Z_{jj'} := {\ts\frac{1}{2}} iZ_{jj'}, \quad
   &\t e_j := {\ts\frac{1}{2}} (f_j + ie_j), \quad
   &\t f_j := {\ts\frac{1}{2}} (f_j - ie_j).
\end{array} \end{equation*}
One checks that there is a unique $\t K^J_N$-splitting $\t\dk^J_N \oplus \t\dm^J_N$ of $\t\dg^J_N$, given by
\begin{equation*}
   \t\dk^J_N  = \Span \bigl\{\t H,\, \t Z_{jj'}:\, 1 \le j \le j' \le N \bigr\}, \quad
   \t\dm^J_N := \Span \bigl\{\t E,\, \t F,\, \t e_j,\, \t f_j:\, 1 \le j \le N \bigr\}.
\end{equation*}
It is simple to verify that $\tau$ is an automorphism, and so the given basis of $\t\dm^J_N$ is a $\t K^J_N$-eigenbasis: the $\t H$-weights of $\t E$, $\t F$, $\t e_j$, and $\t f_j$ are $2$, $-2$, $1$, and $-1$, respectively.

Write $\pi_{k,L}$ for the $\t K^J_N$-character $\pi_{\beta^k \alpha_L}$ associated to the scalar cocycle $\beta^k \alpha_L$ defining the slash action $|_{k,L}$, and $\bC_{k,L}$ for its space.  Check that $\pi_{k,L}$ is determined by $\t H \mapsto -k$ and $\t Z_{jj'} \mapsto \pi L_{jj'}$.  Hence $\t E$, $\t F$, $\t e_j$, and $\t f_j$ span copies of $\bC_{-2,0}$, $\bC_{2,0}$, $\bC_{-1,0}$, and $\bC_{1,0}$, respectively, and so $\t E^{\CDO}$, $\t F^{\CDO}$, $\t e_j^{\CDO}$, and $\t f_j^{\CDO}$ are CDOs from $|_{k,L}$ to $|_{k+2,L}$, $|_{k-2,L}$, $|_{k+1,L}$, and $|_{k-1,L}$, respectively.

We remark that $\tau$ maps the Casimir element $\Omega_N$ to $\bigl( \frac{i}{2} \bigr)^N \Omega_N$.  The second order operators $|_{k,L}[\nu(\t H)]$, $|_{k,L}[\nu(\t E)]$, and $|_{k,L}[\nu(\t F)]$ are the analogs for $N>1$ of the operators $\Delta_1$ and $D_\pm$ given on p.~38 of \cite{BS98}.

\medbreak \noindent {\bf Proof of Proposition~\ref{prop:RLops}.}
By the first paragraph of Corollary~\ref{cor:CDObasis}, $\bigl\{\t E^{\CDO} \bigr\}$, $\bigl\{\t F^{\CDO} \bigr\}$, $\bigl\{\t e_j^{\CDO}:\, 1 \le j \le N \bigr\}$, and $\bigl\{\t f_j^{\CDO}:\, 1 \le j \le N \bigr\}$ are bases of $\bD^1(k, L;\, k+2, L)$, $\bD^1(k, L;\, k-2, L)$, $\bD^1(k, L;\, k+1, L)$, and $\bD^1(k, L;\, k-1, L)$, respectively, and there are no other CDOs of order~1.

By~(\ref{eq:Lie slash action}) and Corollary~\ref{cor:CDObasis}, the symbols of the order~1 CDOs at $(i,0)$ are (using vector notation for $e$ and $f$)
\begin{equation*}
   \t E^{\CDO} \equiv -2i \partial_\tau, \quad
   \t F^{\CDO} \equiv 2i \partial_{\ov \tau}, \quad
   \t e^{\CDO} \equiv -i \partial_z, \quad
   \t f^{\CDO} \equiv i \partial_{\ov z}.
\end{equation*}
Once we prove that $X_\pm$ and $Y_\pm$ really are CDOs, it will follow that
\begin{equation*}
   \t E^{\CDO} = -X_+, \quad
   \t F^{\CDO} = -X_-, \quad
   \t e^{\CDO} = -Y_+, \quad
   \t f^{\CDO} = -Y_-,
\end{equation*}
and the first paragraph of the proposition will be proven.  One could carry this out by applying the map $\CDO_{\gamma, \gamma'}$.  Our method was to guess the CDOs from the $N=1$ case given in \cite{BS98} and check them with a computer.  One can also proceed as follows: note that any order~1 CDO from $|_{k,L}$ to $|_{k',L'}$ lies in the $C^\infty(\HS_{1,N})$-span of $\bigl\{ 1,\, \partial_\tau,\, \partial_{\ov \tau},\, \partial_{z_j}, \partial_{\ov z_j}:\, 1 \le j \le N \bigr\}$, and use covariance to solve for the coefficients.  For example, covariance with respect to $E$ and $e$ implies that the coefficients depend only on $y$ and $v$.  Covariance with respect to $Z$ and $H$ implies that $L' = L$ and the entire operator is of weight $k - k'$ under the Euler operator $2\tau \partial_\tau + 2\ov \tau \partial_{\ov \tau} + z^\T \partial_z + \ov z^\T \partial_{\ov z}$.  The condition for covariance with respect to $F$ and $f$ is harder but can be deduced by hand, leading to the formulas for $X_\pm$ and $Y_\pm$.

To prove the second paragraph of the proposition, use~(\ref{eq:commutators}).  It implies that $X_+$ and $Y_{+,j}$ commute among themselves, $X_-$ and $Y_{-,j}$ commute among themselves, $[X_+, X_-]$ and $[Y_{+,j}, Y_{-,k}]$ are constants, and $[X_\pm, Y_{\mp, j}]$ is in the span of the $Y_{\pm, k}$.  To deduce the constants and the coefficients, apply the commutators to~$1$, or if that fails, to $\ov z$. \hfill $\Box$

\medbreak \noindent {\bf Proof of Proposition~\ref{prop:allCDOs}.}
The first paragraph holds by Corollary~\ref{cor:CDObasis} and the fact that both the raising operators and the lowering operators commute among themselves.  For the second paragraph, use~(\ref{casop}) and match the $\tau$- and $z$-symbols separately to see (easily) that $\cC^{k,L}$ minus the first two lines on the right of~(\ref{eq:casimir_raisinglowering}) is in the span of the operators $Y_{+,j} Y_{-,k}$.  The coefficients can be deduced by applying both sides to $\ov z$.  \hfill $\Box$

\medbreak \noindent {\bf Proof of Proposition~\ref{prop:allIDOs}.}
This follows easily from Proposition~\ref{prop:RLops} and Corollary~\ref{cor:CDObasis}.  \hfill $\Box$


%% file: paper.bbl
\providecommand{\bysame}{\leavevmode\hbox to3em{\hrulefill}\thinspace}
\providecommand{\MR}{\relax\ifhmode\unskip\space\fi MR }
\providecommand{\MRhref}[2]{%
  \href{http://www.ams.org/mathscinet-getitem?mr=#1}{#2}
}
\providecommand{\href}[2]{#2}
\begin{thebibliography}{MVW87}

\bibitem[BS98]{BS98}
R.~Berndt and R.~Schmidt, \emph{Elements of the {R}epresentation {T}heory of the {J}acobi group}, Progress in Mathematics, vol.~163, Birkh\"auser Verlag, Basel, 1998.

\bibitem[Bor76]{Bo76}
W.~Borho, \emph{Primitive und vollprimitive {I}deale in {E}inh\"ullenden von $\mathfrak{so}(5, \mathbb{C})$}, J.\ Algebra \textbf{43} (1976), 619--654.

\bibitem[Bri08]{Br08_1}
K.~Bringmann, \emph{On the explicit construction of higher deformations of partition statistics}, Duke Math.\ J.\ \textbf{144} (2008), no.~2, 195--233.

\bibitem[BCR07]{BCR07}
K.~Bringmann, C.~H.\ Conley, and O.~K.\ Richter, \emph{{M}aass-{J}acobi forms over complex quadratic fields}, Math.\ Res.\ Lett.\ \textbf{14} (2007), 137--156.

\bibitem[BCR12]{BCR10}
\bysame, \emph{Jacobi forms over complex quadratic fields via the cubic {C}asimir operators}, Comment.\ Math.\ Helv.\ \textbf{87} (2012), no.~4, 825--859.

\bibitem[BGM09]{BGM09}
K.~Bringmann, F.~Garvan, and K.~Mahlburg, \emph{Partition statistics and quasiharmonic {M}aass forms}, Internat.\ Math.\ Res.\ Notices\ \textbf{2009}, no.~1, 63--97.

\bibitem[BL09]{BL09}
K.~Bringmann and J.~Lovejoy, \emph{Overpartitions and class numbers of binary quadratic forms}, Proc.\ Natl.\ Acad.\ Sci.\ USA \textbf{106} (2009), no.~14, 5513--5516.

\bibitem[BRR12]{BRR12}
K.~Bringmann, M.~Raum, and O.~K.~Richter, \emph{Kohnen’s limit process for real-analytic {S}iegel modular forms}. Adv.\ Math.\ \textbf{231} (2012), no.~2, 1100--1118.

\bibitem[BRR14]{BRRZ11}
\bysame, \emph{Harmonic Maass-Jacobi forms with singularities and a theta-like decomposition}, Trans.\ Amer.\ Math.\ Soc.\ (2014), in press.

\bibitem[BR10]{BR10}
K.~Bringmann and O.~K.\ Richter, \emph{Zagier-type dualities and lifting maps for harmonic {M}aass-{J}acobi forms}, Adv.\ Math.\ \textbf{225} (2010), 2298--2315.

\bibitem[BZ10]{BZ10}
K.~Bringmann and S.~Zwegers, \emph{Rank-crank type {PDE}'s and non-holomorphic {J}acobi forms}, Math.\ {R}es.\ {L}ett. \textbf{17} (2010), no.~4, 589--600.

\bibitem[BF04]{BF04}
J.~H. Bruinier and J.~Funke, \emph{On two geometric theta lifts}, Duke Math.\ J.\ \textbf{125} (2004), no.~1, 45--90.

\bibitem[CSL09]{CL09}
R.~Campoamor-Stursburg and S.~G. Low, \emph{Virtual copies of semisimple {L}ie algebras in enveloping algebras of semidirect products and {C}asimir operators}, J.\ Phys.\ A: Math.\ Theor. \textbf{42} (2009), 065205 (18pp).

\bibitem[DMZ12]{DMZ11}
A.~Dabholkar, S.~Murthy, and D.~Zagier, \emph{Quantum black holes, wall crossing, and mock modular forms}, {\tt arXiv:1208.4074.}

\bibitem[EZ85]{EZ85}
M.~Eichler and D.~Zagier, \emph{The {T}heory of {J}acobi {F}orms}, Birkh\"auser, Boston, 1985.

\bibitem[Fay77]{Fa77}
J.~D.\ Fay, \emph{Fourier coefficients of the resolvent for a {F}uchsian group}, J.\ Reine Angew.\ Math.\ \textbf{293/294} (1977), 143--203.

\bibitem[Hay06]{Ha06}
S.~Hayashida, \emph{Skew-holomorphic {J}acobi forms of higher degree}, Automorphic {F}orms and {Z}eta {F}unctions, World Scientific, Hackensack, 2006, pp.~130-139.

\bibitem[Hel77]{Helg77}
S.~Helgason, \emph{Invariant differential equations on homogeneous manifolds}, Bull.\ {A}mer.\ {M}ath.\ {S}oc.\ \textbf{84} (1977), no.~5, 751--774.

\bibitem[KW94]{KW94}
V.~G.\ Kac and M.~Wakimoto, \emph{Integrable highest weight modules over affine superalgebras and number theory}, Lie {T}heory and {G}eometry, Progr.\ Math., vol.~123, Birkh\"auser, Boston, 1994, pp.~415--456.

\bibitem[KW01]{KW01}
V.~G.\ Kac and M.~Wakimoto, \emph{Integrable highest weight modules over affine superalgebras and {A}ppell's function}, Comm.\ Math.\ Phys.\ \textbf{215} (2001), no.~3, 631--682.

\bibitem[LV80]{LV80}
G.~Lion and M.~Vergne, \emph{The {W}eil {R}epresentation, {M}aslov {I}ndex, and {T}heta {S}eries}, Progress in Mathematics, vol.~6, Birkh\"auser, Boston, 1980.

\bibitem[Maa49]{Ma49}
H.~Maa\ss, \emph{\"{U}ber eine neue {A}rt von nichtanalytischen automorphen {F}unktionen und die {B}estimmung {D}irichletscher {R}eihen durch {F}unktionalgleichungen}, Math.\ Ann.\ \textbf{121} (1949), 141--183.

\bibitem[MO12]{MO10}
A.~Malmendier and K.~Ono, \emph{{${\rm SO}(3)$}-{D}onaldson invariants of\/ {$\Bbb{C}{\rm P}^2$} and mock theta functions}, Geom.\ Topol.\ \textbf{16} (2012), no.~3, 1767--1833.

\bibitem[MVW87]{MVW87}
C.~M{\oe}glin, M.-F. Vign{\'e}ras, and J.-L.\ Waldspurger, \emph{Correspondances de {H}owe sur un {C}orps {$p$}-adique}, Lecture Notes in Mathematics, vol.~1291, Springer-Verlag, Berlin, 1987.

\bibitem[Ono09]{Ono08}
K.~Ono, \emph{Unearthing the visions of a master: harmonic {M}aass forms and number theory}, Current {D}evelopments in {M}athematics, 2008, Int.\ {P}ress, {S}omerville, 2009, pp.~347--454.

\bibitem[Pit09]{Pi09}
A.~Pitale, \emph{Jacobi {M}aa\ss\ forms}, Abh.\ Math.\ Semin.\ Univ.\ Hambg.\ \textbf{79} (2009), no.~1, 87--111.

\bibitem[Que88]{Qu88}
C.~Quesne, \emph{Casimir operators of semidirect sum {L}ie algebras}, J.\ Phys.\ A: Math.\ Gen. \textbf{21} (1988), L321--L324.

\bibitem[STT05]{STT05}
A.~M. Semikhatov, A.~Taormina, and I.~Yu. Tipunin, \emph{Higher-level {A}ppell functions, modular transformations, and characters}, Comm.\ Math.\ Phys.\ \textbf{255} (2005), no.~2, 469--512.

\bibitem[Sko90]{Sk90}
N.-P. Skoruppa, \emph{Explicit formulas for the {F}ourier coefficients of {J}acobi and elliptic modular forms}, Invent.\ Math.\ \textbf{102} (1990), no.~3, 501--520.

\bibitem[WW27]{WW96}
E.~T. Whittaker and G.~N. Watson, \emph{A {C}ourse of {M}odern {A}nalysis}, $4^\thup$ ed., Cambridge University Press, Cambridge, 1927.

\bibitem[Zwe02]{Zw02}
S.~Zwegers, \emph{Mock theta functions}, Ph.D.\ Thesis, Universiteit Utrecht, 2002.

\bibitem[Zwe10]{Zw10}
\bysame, \emph{Multivariable {A}ppell functions}, Unpublished, 2010.

\end{thebibliography}
